\documentclass[10pt,twocolumn,twoside]{IEEEtran}

\usepackage{graphicx}
\usepackage{cite}
\usepackage{amsmath}
\usepackage{amssymb,amsthm}
\usepackage{color}
\usepackage[figuresright]{rotating}
\usepackage{subfigure}
\usepackage{graphics}
\usepackage[hidelinks]{hyperref}
\usepackage{fontenc}

\usepackage{ifpdf}
\ifpdf
\fi

\newcommand{\tr}[1]{\textbf{tr} \left( #1 \right) }

\newcommand{\Hso}{\ensuremath{H_{\mathrm{SO}}}}
\newcommand{\Hfo}{\ensuremath{H_{\mathrm{FO}}}}

\newcommand{\E}{{\cal E}}
\newcommand{\V}{{\cal V}}

\newcommand{\G}{{\cal G}}
\newcommand{\Z}{\ensuremath{\mathbb{Z}}}

\newtheorem{theorem}{Theorem}[section]
\newtheorem{definition}{Definition}[section]

\newtheorem{proposition}[theorem]{Proposition}

\def\Nsf2{{N}_{\mbox{\scriptsize \rm{2SF}}}}

\begin{document}
\title{Biharmonic Distance  and the Performance of Second-Order Consensus Networks  with Stochastic Disturbances}
 \author{Yuhao~Yi, Bingjia~Yang, Zhongzhi~Zhang, and Stacy~Patterson, \emph{Member, IEEE}
 \thanks{Yuhao Yi is with the Shanghai Key Laboratory of Intelligent Information Processing, School of Computer Science, Fudan University, Shanghai, 200433, China.}
 \thanks{Bingjia Yang is with the Department of Physics, Fudan University, Shanghai, 200433, China, and the Shanghai Key Laboratory of Intelligent Information Processing, School of Computer Science, Fudan University, Shanghai, 200433, China.}
 \thanks{Zhongzhi Zhang is with the Shanghai Key Laboratory of Intelligent Information Processing, School of Computer Science, Fudan University, Shanghai, 200433, China. {\tt\small zhangzz@fudan.edu.cn} }
 \thanks{Stacy Patterson is with the Department of Computer Science, Rensselaer Polytechnic Institute, Troy, New York, 12180.
 {\tt\small sep@cs.rpi.edu}
 }
 }
\maketitle

\begin{abstract}
We study  second order consensus dynamics with random additive disturbances. We investigate three different performance measures:   the steady-state variance of pairwise differences between vertex states, the steady-state variance of the deviation of each vertex state from the average, and the total steady-state variance of the system. We show that these performance measures are closely related to the biharmonic distance;  the square of the biharmonic distance plays similar role in the system performance as  resistance distances  plays in the performance of first-order noisy consensus dynamics. We further define the new concepts of biharmonic Kirchhoff index and vertex centrality based on the biharmonic distance. Finally, we derive analytical results for the performance measures and concepts for complete graphs, star graphs, cycles, and paths, and we use this analysis to  compare the asymptotic behavior of the steady-variance in first- and second-order systems.  
\end{abstract}


\begin{IEEEkeywords}
Distributed average consensus, network coherence, Laplacian spectral distance, biharmonic distances, Gaussian white noise
\end{IEEEkeywords}

%
\IEEEpeerreviewmaketitle

\section{Introduction}
Consensus dynamics have been studied intensively in the context of distributed networked systems because these dynamics represent a fundamental way of sharing information between agents in the network. Consensus algorithms can be widely applied to many real-world applications such as clock synchronization~\cite{CaZa14,SuStBrGh15}, load balancing~\cite{DiFrMo99}, sensor networks~\cite{LiRu06}, formation control~\cite{FaMu04} and distributed optimization~\cite{Sa14}.

In consensus dynamics, when nodes are subject to external disturbances, these disturbances prevent the system from reaching consensus, instead making node states fluctuate around the current average~\cite{BaJoMiPa12}.  Many works have explored analytical methods to quantify the steady-state variance of the deviations from the average. The vast majority of these have considered first-order consensus algorithms~\cite{BaJoMiPa12,YoScLe10,PaBa14,FiLe16,YiZhLiCh15,JaOl15}.  It has been shown that, in such systems,  the total steady-state variance  can be described by  resistance distances in an associated electrical network~\cite{YoScLe10,PaBa14}.  And, in turn, resistance distances are given by the covariance matrix of the vertex states in such a dynamical system~\cite{YoScLe16}.



Many real world systems can be more accurately modeled using second-order dynamics.   For example, second-order consensus protocols are applied to formation control because they capture the kinematics of the vehicles~\cite{ReAt05}. Clock synchronization algorithms using second-order consensus scheme have also been studied~\cite{CaZa14}. 
While second-order dynamics have important applications, analysis of the effects of external perturbations on second-order systems remains limited when compared to recent work on first-order systems.
Previous works have shown that the \emph{total steady-state variance} in such systems are determined by the eigenvalues of the Laplacian matrix, and  asymptotic behaviors for macroscopic and microscopic behaviors of the variance have been so studied in~\cite{BaJoMiPa12}. 
However, no unified metric for second-order systems that is similar to resistance distance for first-order systems has been previously proposed.


In this paper, we propose biharmonic distance as a tool to analyze second-order consensus dynamics with external perturbations.
Biharmonic distance is defined based on the spectrum of the Laplacian matrix, and it has been used in computer graphics~\cite{LiRuFu10}  as a metric that incorporates both local and global graph structure. 
We study three performance measures in second-order consensus systems: the variance of of the difference between the states of any pair of vertices, the variance between an individual vertex state and the system average, and the total variance of the system. 
For each of these performance measures, we show how it can be analyzed in terms of biharmonic distances. 
  In addition, we introduce a new notion of vertex centrality based on a biharmonic vertex index. A vertex with higher biharmonic centrality has smaller steady-state variance.
  We then derive closed-form solutions for the biharmonic distances and related performance measures for complete graphs, star graphs, cycles, and paths.
 Finally,  we use this analysis to  compare the  behavior of the steady-variance in first- and second-order systems.

\subsubsection*{Related work} Bamieh et al. introduced the concept of \emph{network coherence}, a measure of the average steady-state variance of node states, for both first- and second-order consensus dynamics with stochastic external perturbations. This work showed a relationship between coherence and the spectrum of the Laplacian matrix and derived the asymptotic behavior of coherence in torus networks~\cite{BaJoMiPa12}. 
Several works have analyzed the coherence of first-order consensus in different classes of networks. 
Young et al.~\cite{YoScLe10} elated network coherence to the Kirchhoff index of a graph and presented closed-form results for the coherence of cycle, path, and star graphs with first-order noisy consensus dynamics. Patterson and Bamieh analyzed coherence in several forms of fractal trees~\cite{PaBa14} and discussed the impact of fractal dimensions on network coherence, and Yi et al. investigated coherence in Farey graphs~\cite{YiZhLiCh15} and Koch graphs~\cite{YiZhShCh17} as deterministic generated representatives of small-world networks and scale-free networks. 



There have also been several recent works on analysis of coherence for second-order systems in different graph topologies. Namely, the second-order coherence of torus~\cite{BaJoMiPa12}, fractals~\cite{PaBa14}, and Koch graphs~\cite{YiZhShCh17} have all been analyzed. However, none of these works have developed a general mathemtical connection between second-order coherence and a graph distance metric.

With respect to biharmonic distance, the recent work by Fitch and Leonard~\cite{FiLe16} used a slightly different definition of this distance  to describe the centrality of multiple leaders in first-order consensus systems with leader nodes. We show that, while related, this different definition cannot be extended to describe coherence in leader-free second-order consensus networks.

The remainder of this paper is organized as follows. In Section~\ref{pre:sec}, we introduce notation and the system dynamics studied in this paper. In Section~\ref{bih:sec}, we first describe the notion of biharmonic distance and its definition. We then introduce graph indices and vertex centrality based on biharmonic distance. In Section~\ref{second:sec}, we show that biharmonic distance plays a important role in perturbed second-order  consensus dynamics, and we give relationships between coherence performance measures and the biharmonic distance and its derived indices.
In Section~\ref{first:sec}, we compare the relationships between first-order noisy consensus dynamics and resistance distance and second-order noisy consensus dynamics and biharmonic distance.Section~\ref{example:sec} gives closed-form solutions for the coherence performance measures for complete graphs, star graphs, cycles, and paths. In Section \ref{num:sec}, we further investigate these performance measures using numerical examples. Finally, we conclude the paper in Section~\ref{conclusion:sec}.

\section{Preliminaries}\label{pre:sec}
\subsection{Concepts and Notation}

Let $\G$ be an undirected connected graph, and let $\V=\{0,1\dots,N-1\}$ and $\E$ be the vertex set and edge set that constitute $\G$ as $\G=\{\V,\E\}$. Let $N = |\V|$ and $M=|\E|$. Define $A$ as the $N\times N$ ($0$-indexed) adjacency matrix of $\G$, in which $a_{ij}=1$ if $\{i,j\}\in E$ and $a_{ij}=0$ otherwise. Let $D$ be the diagonal matrix where $d_{ii}$ is equal to the degree of vertex $i$, \emph{i.e.}, $d_{ii}=\sum_{i=0}^{N-1}a_{ij}$. Define $L=D-A$ as the Laplacian matrix of graph $G$. We use $\lambda_i$ and $u_i$ to denote the $i$-th eigenvalue and eigenvector of $L$, $i\in \{0,1,\dots, N-1\}$, where $0=\lambda_0<
\lambda_1\leq \dots \leq \lambda_{N-1}$. The all-one vector of order $N$ is denoted by $\mathbf{1}_N$.  Therefore, $u_0=\frac{1}{\sqrt{N}}\textbf{1}_N$. Then, $L$ can be diagonalized as $L=U\Lambda U^{\top}$, where $\Lambda\in \mathbb{R}^{N\times N}$ is diagonal and $\Lambda_{ii}=\lambda_i$, $U\in \mathbb{R}^{N\times N}$, with its $i$th column being $u_i$. In addition, we denote by $L^\dagger$ the pseudo-inverse of $L$, and define $L^{2\dagger}=(L^\dagger)^2$.

\subsection{System Dynamics}
Each vertex in the network has a scalar-valued state.
Let $x_1(t)$ be the  $N$-vector that contains the states of all vertices; $x_{1j}(t)$ represents the state of vertex $j$, $j\in\{0,1\dots,N-1\}$. Then, we define $x_2(t)$ as the first derivative of $x_1(t)$ with respect to $t$, that is, $x_2(t)=\dot{x}_1(t)$. A vertex $j$ adjust its state by setting $\dot{x}_{2j}(t)$ according to the differences of its state ($x_{1j}(t)$ and $x_{2j}(t)$) and the states of its neighbors. The following equation gives the noisy second-order consensus algorithm:
\begin{equation} \label{system:eqn}
\left[
\begin{array}{c}
\dot{x}_1(t)\\
\dot{x}_2(t)
\end{array}
\right]
=\left[
\begin{array}{cc}
0 & I\\
-L & -L
\end{array}
\right]
\left[
\begin{array}{c}
x_1(t)\\
x_2(t)
\end{array}
\right]
+\left[
\begin{array}{c}
0\\
I
\end{array}
\right]w(t)\,,
\end{equation}
where $0$, $I$, and $L$ are all $N\times N$ matrices, and $w(t)$ is a $2N$-vector of uncorrelated Gaussian white noise processes.


\subsection{Performance Measures} \label{perfmeas:sec}
Because the state of each vertex is disturbed by  Gaussian noise, the networked system can never reach exact consensus. Therefore, we are interested in the expected deviations of the states of the vertices. In particular, we are interested in three performance measures related to these deviations, which we define below.

First, we want to know how far the states of two vertices are driven away by disturbances. Therefore we study the steady-state of the variance of this pairwise deviation.
\begin{definition}
For any two vertices $j, k \in \V$, the \emph{pairwise variance} $\Hso(j,k)$ is the steady-state variance of the difference between $x_{1j}$ and $x_{1k}$, i.e.,
\begin{align}
\label{hso_pairwise:def}
\Hso(j,k)=\lim_{t\to\infty}\mathbb{E}[\left(x_{1j}(t)-x_{1k}(t)\right)^2].
\end{align}
\end{definition}
We note that in a $d$-dimensional torus $\Z^d_N$, $\Hso(j,j-1)$ is the second-order microscopic coherence defined in~\cite{BaJoMiPa12}, and $\Hso(j,j+\frac{N}{2})$ is the second-order long-range coherence defined in~\cite{BaJoMiPa12}.
Thus, our pairwise variance performance measure is a generalization of these two performance measures.

We are also interested in the variance of the difference between the state of a vertex and the (current) average value in the network. Let $\bar{x}_1(t)$ be the average state $\bar{x}_1(t)=\frac{1}{N} \mathbf{1}_N^{\top}x_{1}(t)$.
\begin{definition}
For a vertex $j \in \V$, the \emph{vertex variance} $\Hso(j)$ is the steady-state variance of the difference between $x_{1j}(t)$ and $\bar{x}_1(t)$, i.e.,\begin{align}
\label{hso_vertex:def}
\Hso(j)=\lim_{t\to\infty}\mathbb{E}[\left(x_{1j}(t)-\bar{x}_{1}(t)\right)^2]\,.
\end{align}
\end{definition}

Finally, we are  interested in the total variance of the system. 
\begin{definition}
For a network $\G$, the \emph{total variance} $\Hso(\G)$ is the total steady-state variance of the deviation of each vertex state from the current average, i.e., 
\begin{align}
\label{hso_total:def}
\Hso(\G)=\lim_{t\to\infty}\sum_{j=0}^{N-1}\mathbb{E}[\left(x_{1j}(t)-\bar{x}_{1}(t)\right)^2]\,.
\end{align}
\end{definition}
In a $d$-dimensional torus $\Z^d_N$, $\Hso(\G)$ is the variance of the deviation from average defined in~\cite{BaJoMiPa12}.

\section{Biharmonic Distance}\label{bih:sec}

Several slightly different definitions of biharmonic distance have been proposed in related literature~\cite{LiRuFu10,Pa17book,FiLe16}. In this paper we follow the definition in~\cite{LiRuFu10} and~\cite{Pa17book}, which is as follows.
\begin{definition}
\label{biharmonic:def}
The biharmonic distance $d_B(j,k)$ between two vertices $j$ and $k$ in a undirected graph $\G$ is:
\begin{small}
\begin{align}
\label{biharmonic_def:eqn}
d^2_B(j,k)=L^{2\dagger}_{jj}+L^{2\dagger}_{kk}-2L^{2\dagger}_{jk}=\sum_{i=1}^{N-1}\frac{1}{\lambda_i^2}(u_{ij}-u_{ik})^2\,.
\end{align}
\end{small}
\end{definition}
Note that this definition is equal to the square root of the one used by Fitch and Leonard in~\cite{FiLe16}.

Biharmonic distance is a metric, as shown in the following theorem. While this result has been previously proved~\cite{LiRuFu10}, we include a proof for the convenience of the reader.
\begin{theorem}
The biharmonic distance $d_B(j,k)$ is a $\V\times \V\to \mathbb{R}$ metric, which is equivalent to satisfying the following properties:
\begin{itemize}
    \item Non-negativity: $d_B(j,k)\geqslant 0 $,
    \item Nullity: $d_B(j,k)=0$ if and only if $j=k$,
    \item Symmetry: $d_B(j,k)=d_B(k,j)$, and
    \item Triangle inequality $d_B(j,r)+d_B(r,k)\geqslant d_B(j,k)$.
\end{itemize}
\end{theorem}
\begin{IEEEproof}
The non-negativity and symmetry are easily obtained from Definition~\ref{biharmonic:def} along with the fact that $L$ is positive semi-definite. Assume $d_B(j,k)=0$ for $j\neq k$, then $u_{ij}=u_{ik}$ for all $i\in\{0,1,\dots,N-1\}$. Since $L=U\Lambda U^{\top}$,  $L_{jj}=\sum_{i=1}^{N-1}\lambda_i u_{ij}u_{ij}$ and $L_{jk}=\sum_{i=1}^{N-1}\lambda_i u_{ij}u_{ik}$. This leads to $L_{jj}=L_{jk}$ for $j\neq k$, which contradicts with the definition of the Laplacian matrix. 

The triangle inequality can be proved as follows. Define a vector,
\begin{align}
v_j=\sum_{i=1}^{N-1}\frac{u_{ij}}{\lambda_i}u_i \in \mathbb{R}^N \text{ for } j=0, 1,\dots,N-1\,.\nonumber
\end{align}
We note again that $u_i$ is the $i$th eigenvector, and $u_{ij}$ is the $j$th entry of $u_i$. Then it follows that 
the Euclidean distance $\|v_j-v_k\|_2$ between $v_j$ and $v_k$ is
$$\|v_j-v_k\|_2=\left\|\sum_{i=1}^{N-1}\frac{(u_{ij}-u_{ik})}{\lambda_i}u_i\right\|_2=\sqrt{\sum_{i=1}^{N-1}\frac{(u_{ij}-u_{ik})^2}{\lambda^2_i}}\,,$$
which means $d_B(j,k)$ is equal to $\|v_j-v_k\|_2$. Since the Euclidean distance in $\mathbb{R}^N$ is a metric and, therefore, satisfies triangle inequality, $d_B(j,k)$ also satisfies the triangle inequality.
\end{IEEEproof}
We observe that $v_j, j\in\{1,\dots,N\}$ assigns a position to vertex $j$ in $\mathbb{R}^N$ Euclidean space that preserves biharmonic distance.
\begin{definition}
We define an $N$-dimensional mapping of of the vertices in $\G$, $\mathcal{F}: \V\to \mathbb{R}^N$. For any vertex $j$, ${\mathcal{F}}(j)=v_j=L^{\dagger}e_j$. $v_j$ is a \emph{biharmonic embeddings} of graph $\G$ in $\mathbb{R}^N$.
\end{definition}

Based on the definition of biharmonic distance, we  also define the following graph indices.
\begin{definition}
The \emph{biharmonic Kirchhoff index} $D^2_{B}(\G)$ of a graph $\G$ is 
\begin{align}
D^2_{B}(\G) = \sum_{\substack{j,k\in V\\j<k}}d^2_B(j,k)\,.
\end{align}
\end{definition}
\begin{definition}
The \emph{biharmonic vertex index} $D^2_{B}(j)$ of a node $j$ in a graph $\G$ is 
\begin{align}
D^2_B(j)=\sum_{k\in V}d_B^2(j,k)\,.
\end{align}
\end{definition}
We can derive from the definition of $d_B(j,k)$ that
\begin{align}
\label{dbg_eigen:eqn}
D^2_{B}(\G)=N\cdot \sum_{i=1}^{N-1}\frac{1}{(\lambda_i)^2}\,.
\end{align}


Finally, for a vertex $j$ in graph $\G$, we can define its centrality based on  biharmonic distances.
\begin{definition}
The \emph{biharmonic centrality} of vertex $j$ in graph $\G$ is
\begin{align}
C_B(j)=\left(\frac{1}{N}D^2_B(j)\right)^{-1}\,.
\end{align}
\end{definition}


\section{Biharmonic Distance in Second-order Consensus Dynamics with Disturbances}\label{second:sec}

The equation
(\ref{system:eqn}) gives  the dynamics of the second-order consensus algorithm with stochastic perturbations. 
The deviation of the state of vertex $j$ from the average of all states is given by $y_j(t)=x_{1j}(t)-\bar{x}_1(t)$. Let $y(t)$ be a $N\times 1$ vector representing all vertices' deviations from average,
\begin{align*}
y(t)=\left[\Pi~|~0 \right] x(t)=\Pi x_1(t)\,,
\end{align*}
where $\Pi=I_{N}- \frac{1}{N}\mathbf{1}_N \mathbf{1}^\top_N$. The performance measures we study in this paper can all be expressed in terms of  of $y(t)$. Specifically, 
\begin{small}
\begin{align}
\Hso(j,k)=&\lim_{t\to\infty}\mathbb{E}[\left((x_{1j}(t)-\bar{x}_1(t))-(x_{1k}(t)-\bar{x}_1(t))\right)^2]\nonumber\\
=&\lim_{t\to\infty}\mathbb{E}[\left(y_{j}(t)-y_{k}(t)\right)^2]\\
\Hso(j)=&\lim_{t\to\infty}\mathbb{E}[\left(x_{1j}(t)-\bar{x}_{1}(t)\right)^2]=\lim_{t\to\infty}\mathbb{E}[y_j(t)^2]\\
\Hso(\G)=&\lim_{t\to\infty}\sum_{j=0}^{N-1}\mathbb{E}[y_j(t)^2]\,.
\end{align}
\end{small}

However, the system described by (\ref{system:eqn}) is only marginally stable~\cite{YoScLe10}.  To obtain a stable system, 
we only consider the dynamics in the subspace that is orthogonal to the subspace spanned by $\mathbf{1}_N$. 
We define $Q$ as a $(N-1)\times N$ 
matrix whose rows are the eigenvectors of $L$, excluding $\mathbf{1}_N$. 
We recall that $L$ can be diagonalized as $U\Lambda U^{\top}$, where $U$ is a unitary matrix and $\Lambda$ is a diagonal matrix. Then, $Q^\top$ is the submatrix of $U$ formed by eliminating the first column. It is easy to confirm that $Q\mathbf{1}_N=0$, $QQ^\top=I_{N-1}$, $Q^{\top}Q=\Pi$, and $LQ^{\top}Q=L$. Then, we define
\begin{align*}
z_1(t)=\left[~Q~|~0~\right] x(t)=Q x_1(t)\,,
\end{align*}
and note that $y(t)=Q^{\top}z_1(t)$. It indicates that we can write expressions for our performance measures using $z_1(t)$. Let $z_2(t)=\dot{z}_1(t)$. 
Then (\ref{system:eqn}) leads to
\begin{small}
\begin{align*}
\left[
\begin{array}{cc}
Q & 0\\
0 & Q
\end{array}
\right]
\left[
\begin{array}{c}
\dot{x}_1(t)\\
\dot{x}_2(t)
\end{array}
\right]
=&\left[
\begin{array}{cc}
0 & Q\\
-QLQ^{\top}Q & -QLQ^{\top}Q
\end{array}
\right]
\left[
\begin{array}{c}
x_1(t)\\
x_2(t)
\end{array}
\right]\nonumber\\
&+\left[
\begin{array}{c}
0\\
Q
\end{array}
\right]w(t)\,,
\end{align*}
\end{small}
Therefore, we obtain a stable system:
\begin{align*}
\left[
\begin{array}{c}
\dot{z}_1(t)\\
\dot{z}_2(t)
\end{array}
\right]
=\left[
\begin{array}{cc}
0 & I_{N-1}\\
-\bar{\Lambda} & -\bar{\Lambda}
\end{array}
\right]
\left[
\begin{array}{c}
z_1(t)\\
z_2(t)
\end{array}
\right]
+\left[
\begin{array}{c}
0\\
Q
\end{array}
\right]w(t)\,,
\end{align*}
where $\bar{\Lambda}=QLQ^{\top}=QU\Lambda (QU)^{\top}=diag(\lambda_1,\dots,\lambda_{N-1})$.
%

We can always find the unitary (orthogonal) permutation matrix $V\in\{0,1\}^{(2N-2)\times (2N-2)}$ such that
\begin{align}
\label{permut:def}
\left[
\begin{array}{c}
\dot{z}_1(t)\\
\dot{z}_2(t)
\end{array}
\right]
=&
V^{\top} K V \left[
\begin{array}{c}
z_1(t)\\
z_2(t)
\end{array}
\right]
+\left[
\begin{array}{c}
0\\
Q
\end{array}
\right]w(t)\,,
\end{align}
where $K$ is the block diagonal matrix,
\begin{align}
\label{kblockDiag:def}
K = \left[
\begin{array}{c|c|c}
P_1 & &\\
\hline
&\ddots &\\
 \hline
 & & P_{N-1}
\end{array}
\right],
\end{align}
with each $P_i$ defined as: 
\begin{align*}
P_i=\left[
\begin{array}{cc}
0 & 1\\
-\lambda_i &- \lambda_i
\end{array}
\right]\,.
\end{align*}
Hereafter, we use the system dynamics in (\ref{permut:def}) to develop expressions for the performance measures defined in Section~\ref{perfmeas:sec}.

\subsection{Pairwise Variance}

\begin{theorem}
The pairwise variance of the difference between states of vertices $j$ and $k$ with dynamics~(\ref{system:eqn}) can be expressed by the spectrum of the Laplacian matrix of graph $\G$ as
\begin{align}
\Hso(j,k)=\sum_{i=1}^{N-1}\frac{(u_{ij}-u_{ik})^2}{2\lambda_i^2}\,.
\end{align}
\end{theorem}
\begin{IEEEproof}
We start by expressing $\Hso(j,k)$ in terms of $z_1(t)$,
\begin{align*}
&\Hso(j,k)
=\lim_{t\to \infty}\mathbb{E}\left[y(t)^{\top} (e_j-e_k) (e_j^{\top}-e_k^{\top}) y(t)\right]\\
&=\lim_{t\to \infty}\mathbb{E}\left[(Q^{\top} z_1(t))^{\top} (e_j-e_k) (e_j^{\top}-e_k^{\top}) Q^{\top}z_1(t)\right]\\
&=\lim_{t\to \infty}\mathbb{E}\left[z_1(t)^{\top} Q (e_j-e_k) (e_j^{\top}-e_k^{\top}) Q^{\top} z_1(t)\right]\\
&=\lim_{t\to \infty}\mathbb{E}\left[\mathbf{tr}((e_j-e_k)^{\top}Q^{\top}z_1(t) z_1(t)^{\top}Q (e_j-e_k))\right]\,,
\end{align*}
where $e_j$ is the $j$th canonical  basis vector of $\mathbb{R}^{N}$. We define the output of the system as
\begin{align}
\label{pairwiseOutput:def}
\phi(t)&=(e_j-e_k)^{\top}Q^{\top}[I_{N-1} | 0_{N-1}]z(t)\nonumber\\
&= (e_j-e_k)^{\top}Q^{\top} z_1(t)\,.
\end{align}
Then, we define $\Sigma(t)=\mathbb{E}[\phi(t)\phi(t)^{\top}]$; therefore, $\Hso(j,k)=\lim_{t\to\infty}[\tr{\Sigma(t)}] = [\tr{\lim_{t\to\infty}\Sigma(t)}]=:[\tr{\Sigma}]$.

For the state-space system given by~(\ref{permut:def}) and~(\ref{pairwiseOutput:def}), the square of the $\mathcal{H}_2$ norm of the system is
\begin{align}
\label{h2norm:def}
\mathcal{H}^2_2=\int_{0}^{\infty}B^{\top}\mathrm{e}^{-M^{\top}t}Z\mathrm{e}^{-Mt}B \mathrm{d}t\,,
\end{align}
in which
\begin{align}
\label{bexp:eqn}
B=&\left[\begin{array}{c}
0\\
Q
\end{array}\right]\\
\label{mexp:eqn}
M=&\left[
\begin{array}{cc}
0 & I\\
-\bar{\Lambda} & -\bar{\Lambda}
\end{array}
\right]\text{ and }\\
Z=&\left[
\begin{array}{cc}
Q(e_j-e_k) (Q (e_j-e_k))^{\top} & 0\\
0 & 0
\end{array}
\right]\,.
\end{align}
It follows that $\Hso(j,k)=\mathcal{H}^2_2=\mathbf{tr}\left(B^\top \Sigma B\right)$. $\Sigma$ is the solution of the following Lyapunov equation,
\begin{align}
\label{pairLyap:eqn}
M^{\top}\Sigma+\Sigma M+Z=0\,.
\end{align}
The equation is equivalent to
\begin{align*}
&V M^{\top}\Sigma V^{\top}+V\Sigma M V^{\top}=-VZV^{\top}\text{ or }\\
&(V M^{\top}V^{\top})( V\Sigma V^{\top})+(V\Sigma V^{\top})(V M V^{\top})=-VZV^{\top}
\end{align*}
where $V$ was defined in (\ref{permut:def}) as a (unitary) permutation matrix. We denote by $K=V M V^{\top}$ and $\Theta= V\Sigma V^{\top}$. Then equation~(\ref{pairLyap:eqn}) can be written as
\begin{align}
\label{pairModLyap:eqn}
&K^{\top}\Theta+\Theta K=-V Z  V^{\top}\nonumber\\
&=-\left[
\begin{array}{c|c|c}
Z_{11} & \cdots & Z_{1(N-1)} \\
 \hline
\vdots & \ddots & \vdots \\
 \hline
Z_{(N-1)1} &\cdots & Z_{(N-1)(N-1)}
\end{array}
\right]\,,
\end{align}
for $i,m\in \{1,\dots ,N-1\}$,
\begin{align}
& Z_{im}=\left[
\begin{array}{cc}
(Q_{ij}Q_{mj}-Q_{ij}Q_{mk}-Q_{ik}Q_{mj}+Q_{ik}Q_{mk}) & 0\\
0 & 0
\end{array}
\right]\nonumber\\
&=\left[
\begin{array}{cc}
(u_{ij}u_{mj}-u_{ij}u_{mk}-u_{ik}u_{mj}+u_{ik}u_{mk}) & 0\\
0 & 0
\end{array}
\right]\nonumber
\end{align}
We note that $K$ is block-diagonal.
Substituting (\ref{kblockDiag:def}) into diagonal blocks of~(\ref{pairModLyap:eqn}) yields $P^{\top}_i\Theta_{ii}+\Theta_{ii}P_i=Z_{ii}$. Since $Z_{ii}$ and $P_i$ are symmetric, $\Theta_{ii}$ is also symmetric. We write $\Theta_{ii}$ as
\begin{align*}
\Theta_{ii}=\left[\begin{array}{cc}
X_{ii} & \Psi_{ii}\\
\Psi_{ii} & Y_{ii}
\end{array}\right]\,.
\end{align*}
Then,
\begin{align*}
\left[\begin{array}{cc}
0 & \lambda_i\\
1 & \lambda_i
\end{array}\right]
\left[\begin{array}{cc}
X_{ii} & \Psi_{ii}\\
\Psi_{ii} & Y_{ii}
\end{array}\right]
+&\left[\begin{array}{cc}
X_{ii} & \Psi_{ii}\\
\Psi_{ii} & Y_{ii}
\end{array}\right]
\left[\begin{array}{cc}
0 & 1\\
\lambda_i & \lambda_i
\end{array}\right]\nonumber\\
&=\left[\begin{array}{cc}
(u_{ij}-u_{ik})^2 & 0\\
0 & 0
\end{array}\right]\,,
\end{align*}
which leads to
\begin{align*}
Y_{ii}=\frac{(u_{ij}-u_{ik})^2}{2\lambda^2_i}\,.
\end{align*}
Then, we derive that
\begin{align}
\Hso(j,k)=\mathcal{H}^2_2=&\mathbf{tr}\left(B^\top \Sigma B\right)^{\frac{1}{2}}=\mathbf{tr}\left(B^\top V^\top\Theta V B\right)\nonumber\\
&=\sum_{i=1}^{N-1}(Y_{ii})=\sum_{i=1}^{N-1}\frac{(u_{ij}-u_{ik})^2}{2\lambda^2_i}\,.
\end{align}

\end{IEEEproof}

Applying (\ref{biharmonic_def:eqn}), we immediately obtain the following theorem.
\begin{theorem}
\label{hsojk_bih:them}
For any vertex pair $j$ and $k$ in a network $\G$ with dynamics~(\ref{system:eqn}),
\begin{align}
\Hso(j,k)=\frac{1}{2}d^2_{B}(j,k)
\end{align}
\end{theorem}
This theorem shows that the pairwise variance between vertices $j$ and $k$ is proportional to the square of their biharmonic distance.

\subsection{Vertex Variance}
We first give an expression for the vertex variance in terms of the eigenvalues and eigenvectors of $L$.
\begin{theorem} \label{vertexspectrum.thm}
For any vertex $j$ in network $G$ with dynamics~(\ref{system:eqn})
\begin{align}
\label{hsoj_eigen:eqn}
\Hso(j)=\sum_{i=1}^{N-1}\frac{u^2_{ij}}{2\lambda^2_i}\,.
\end{align}
\end{theorem}
\begin{IEEEproof}
First, we derive an expression for the vertex variance in terms of $z_1(t)$,
\begin{align*}
\Hso(j)
&=\lim_{t\to \infty}\mathbb{E}\left[y(t)^{\top} e_j e_j^{\top} y(t)\right]\\
&=\lim_{t\to \infty}\mathbb{E}\left[(Q^{\top} z_1(t))^{\top} e_j e_j^{\top} Q^{\top}z_1(t)\right]\\
&=\lim_{t\to \infty}\mathbb{E}\left[z_1(t)^{\top} Q e_j e_j^{\top} Q^{\top} z_1(t)\right]\\
&=\lim_{t\to \infty}\mathbb{E}\left[\mathbf{tr}(e_j^{\top}Q^{\top}z_1(t) z_1(t)^{\top}Q e_j)\right]\,.
\end{align*}
With this, we define the output for the dynamics (\ref{permut:def}) as,
\begin{align}
\label{vertexOutput:def}
\phi(t)=e_j^{\top}Q^{\top}[I_{N-1} | 0_{N-1}]z(t)= e_j^{\top}Q^{\top} z_1(t)\,.
\end{align}

Again, we define $\Sigma(t)=\mathbb{E}[\phi(t)\phi(t)^{\top}]$, therefore $\Hso(j)=\lim_{t\to\infty}[\tr{\Sigma(t)}] = [\tr{\lim_{t\to\infty}\Sigma(t)}]=:[\tr{\Sigma}]$.

For the state-space system given by~(\ref{permut:def}) and~(\ref{vertexOutput:def}), the square of $\mathcal{H}_2$ norm of the system is also defined by~(\ref{h2norm:def}), in which $B$ and $M$ are given by (\ref{bexp:eqn}) and (\ref{mexp:eqn}), $Z$ is expressed by
\begin{align*}
Z
=\left[
\begin{array}{cc}
Qe_j (Q e_j)^{\top} & 0\\
0 & 0
\end{array}
\right]\,.
\end{align*}
It follows that $\Hso(j)=\mathcal{H}^2_2=\mathbf{tr}\left(B^\top \Sigma B\right)$. $\Sigma$ is the solution of the following Lyapunov equation,
\begin{align}
\label{vertexLyap:eqn}
M^{\top}\Sigma+\Sigma M+Z=0\,,
\end{align}

The equation is equivalent to
\begin{align}
&K^{\top}\Theta+\Theta K=-V Z  V^{\top}\nonumber\\
&\qquad=-\left[
\begin{array}{c|c|c}
Z_{11} & \cdots & Z_{1(N-1)} \\
 \hline
\vdots & \ddots & \vdots \\
 \hline
Z_{(N-1)1} &\cdots & Z_{(N-1)(N-1)}
\end{array}
\right]\,,
\label{Lyapu3}
\end{align}
where
\begin{align*}
Z_{im}=\left[
\begin{array}{cc}
Q_{ij}Q_{mj} & 0\\
0 & 0
\end{array}
\right]=\left[
\begin{array}{cc}
u_{ij}u_{mj} & 0\\
0 & 0
\end{array}
\right]\,,
\end{align*}
for $i,m\in\{1,\dots,N-1\}$. We recall that $K=V M V^{\top}$ and $\Theta=V\Sigma V^{\top}$.

Substituting (\ref{kblockDiag:def}) into diagonal blocks of~(\ref{Lyapu3}) yields $P^{\top}_i\Theta_{ii}+\Theta_{ii}P_i=Z_{ii}$. Similar to the pairwise case, we assume
\begin{align}
\Theta_{ii}=\left[\begin{array}{cc}
X_{ii} & \Psi_{ii}\\
\Psi_{ii} & Y_{ii}
\end{array}\right]\,.
\end{align}
By solving $P^{\top}_i\Theta_{ii}+\Theta_{ii}P_i=Z_{ii}$ we derive
\begin{align*}
Y_{ii}=\frac{u^2_{ij}}{2\lambda^2_i}\,.
\end{align*}
Then we obtain 
\begin{align}
\Hso(j)=\mathcal{H}^2_2=&\mathbf{tr}\left(B^\top \Sigma B\right)=\mathbf{tr}\left(B^\top V^\top\Theta V B\right)\nonumber\\
&=\sum_{i=1}^{N-1}(Y_{ii})=\sum_{i=1}^{N-1}\frac{u^2_{ij}}{2\lambda^2_i}\,.
\end{align}

\end{IEEEproof}

We next use Theorem~\ref{vertexspectrum.thm} to derive an expression for the vertex variance in terms of biharmonic distances.
\begin{theorem}
\label{hsoj_bih:them}
For any vertex $j$ in network $\G$ with dynamics~(\ref{system:eqn}), the variance of difference between the state of a vertex and the system averge is decided by the spectrum of the Laplacian marix of the graph, that is
\begin{align}
\Hso(j)=\frac{1}{2N}\left(D^2_B(j)-\frac{1}{N}D^2_B(\G)\right)\,.
\end{align}
\end{theorem}
\begin{IEEEproof}
The biharmonic distance from vertex $j$ to all other vertices is
\begin{align}
\label{hsoj_eig_expand:eqn}
D^2_B(j)&=\sum_{k=0}^{N-1} d^2_B(j,k)
=\sum_{k=0}^{N-1}\sum_{i=1}^{N-1}\frac{1}{\lambda_i^2}(u_{ij}-u_{ik})^2\nonumber\\
&=\sum_{i=1}^{N-1}\sum_{k=0}^{N-1}\frac{u^2_{ij}-2u_{ij}u_{ik}+u_{ik}^2}{\lambda_i^2}\nonumber\\
&=N\sum_{i=1}^{N-1}\frac{u^2_{ij}}{\lambda_i^2}+\sum_{i=1}^{N-1}\frac{1}{\lambda_i^2}\,.
\end{align}
Substituting (\ref{dbg_eigen:eqn}) and (\ref{hsoj_eigen:eqn}) into (\ref{hsoj_eig_expand:eqn}), we obtain
\begin{align}
\Hso(j)=\frac{D^2_B(j)}{2N}-\frac{D^2_B(\G)}{2N^2}\,.
\end{align}
\end{IEEEproof}

\subsection{Total Variance}
Finally, we present expressions for the total variance in terms of the spectrum of the Laplacian matrix.
\begin{theorem}
The total steady-state variance $\Hso(\G)$ of system (\ref{system:eqn}) is
\begin{equation}
    \Hso(\G)=\sum_{i=0}^{N-1}\frac{1}{2\lambda^2_i}\,.
\end{equation}
\end{theorem}
\begin{IEEEproof}
Since,
\begin{align*}
H_{\text{SO}}(\G)=\sum_{j=0}^{N-1}\Hso(j)\,,
\end{align*}
we immediately obtain
\begin{align*}
\Hso(\G)=\sum_{j=0}^{N-1}\sum_{i=0}^{N-1}\frac{u_{ij}^2}{2\lambda_i^2}=\sum_{i=0}^{N-1}\sum_{j=0}^{N-1}\frac{u_{ij}^2}{2\lambda_i^2}=\sum_{i=0}^{N-1}\frac{1}{2\lambda_i^2}.
\end{align*}
\end{IEEEproof}

In similar fashion, we use (\ref{dbg_eigen:eqn}) to obtain the following theorem about the relationship between the total variance and biharmonic distances.
\begin{theorem}
\label{hsog_bih:them}
For a network $\G$ with dynamics (\ref{system:eqn}), the total  variance is given by the biharmonic Kirchhoff index of the graph, specifically,
\begin{equation}
    \Hso(G)=\frac{1}{2N}D_B^2(\G)\,.
\end{equation}
\end{theorem}

\section{Resistance Distance in First-order Consensus Dynamics with Disturbances} \label{first:sec}
In this section, we briefly review  first-order consensus dynamics with stochastic disturbances
and the relationship between resistance distance and the total steady-state variance

The first-order consensus system is formulated as
\begin{align}
\label{fosystem.eqn}
\dot{x}(t)=-Lx(t)+w(t)\,,
\end{align}\
where $x(t)\in \mathbb{R}^{N}$ represents the states of the vertices, and $w(t)\in \mathbb{R}^{N}$ is a vector of uncorrelated Gaussian white noise processes.
The total steady-state variance of the system is
\begin{align}
\Hfo(\G)=\lim_{t\to\infty}\sum_{j=1}^{N}\mathbb{E}[\left(x_{j}(t)-\bar{x}(t)\right)^2]\,,
\end{align}
where $\bar{x}(t)=\frac{1}{N} \mathbf{1}_N^{\top}x(t)$.


The total steady-state variance $\Hfo$ can be expressed in terms of resistance distances in an electrical network. We first formalize the notion of resistance distance and the Kirchhoff index.
\begin{definition}
The resistance distance $d_R(j,k)$ between two vertices $j$ and $k$ in an undirected graph $\G$ is defined as
\begin{small}
\begin{align}
d_R(j,k)=L^{\dagger}_{jj}+L^{\dagger}_{kk}-2L^{\dagger}_{jk}=\sum_{i=1}^{N-1}\frac{1}{\lambda_i}(u_{ij}-u_{ik})^2\,.
\end{align}
\end{small}
\end{definition}
\begin{definition}
The \emph{Kirchhoff index} $D_{R}(\G)$ of a graph $\G$ is defined as
\begin{align}
D_{R}(\G) = \sum_{\substack{j,k\in V\\j<k}}d_R(j,k)\,.
\end{align}
\end{definition}
It has been shown~\cite{BaJoMiPa12,HuSzKo12} that the Kirchhoff index is related to the total steady-state variance of system (\ref{fosystem.eqn}) as
\begin{align}
\Hfo(\G)=\frac{1}{2N}D_R(\G)\,.
\end{align}

We also note that the notion of the \emph{information centrality} of a vertex can be expressed in terms of resistance distances.
If we  defined the sum of resistance distances between all vertices to a vertex $j$ as
\begin{align}
D_R(j)=\sum_{k\in V}d_R(j,k)\,,
\end{align}
then the \emph{information centrality} of vertex $j$ in graph $\G$ is~\cite{FiLe13}
\begin{align}
C_R(j)=\left(\frac{1}{N}D_R(j)\right)^{-1}\,.
\end{align}

Finally, we define the resistance embedding of a graph.
\begin{definition}
Let $\mathcal{F}_R: \V\to \mathbb{R}^N$ be an $n$-dimensional maping of $\G$, such that
for any vertex $j$, ${\mathcal{F}}_R(j)=\mu_j=L^{\dagger/2}e_j$. $\mu_j$ is a \emph{resistance embedding} of graph $\G$ in $\mathbb{R}^N$.
\end{definition}


\section{Analytical Examples}
\label{example:sec}
In this section we give examples for biharmonic distance, connectivity and centrality in networks with special topology. Closed form expressions are derived for all cases. 
We also compare the asymptotic behavior of the steady-state variance of first- and second-order systems.

We note that in some of these examples eigenvectors, are given as complex vectors (although they can be given as real vectors by an unitary linear transform). Therefore, we calculate the  biharmonic distances using the following expression:
\begin{small}
\begin{align}
d_B^2(j,k)=L^{2\dagger}_{jj}+L^{2\dagger}_{kk}-2L^{2\dagger}_{jk}=\sum_{n=1}^{N-1}\frac{1}{\lambda_n^2}|u_{nj}-u_{nk}|^2\,, \label{biharmimg.eq}
\end{align}
\end{small}
which is a slight variation of the definition in (\ref{biharmonic_def:eqn}).
We note that  $i$ is used to indicate the imaginary unit in this section.

\subsection{Complete Graph}
A complete graph is a network in which every vertex is connected to every other vertex. We consider  a complete graph of $N$ vertices.
Its Laplacian matrix of it is
\begin{equation*}
L^{cp}_N=\left(
\begin{array}{ccccc}
  {N-1}& -1&{\cdots}&-1&-1  \\
  -1&{N-1}&{\cdots}&-1&-1 \\ 
  {\vdots}&{\vdots}&{\ddots}&{\vdots}&{\vdots} \\ 
  -1&-1&{\cdots}&{N-1}&-1 \\ 
  -1&-1&{\cdots}&-1&{N-1}\\  
\end{array}
\right) .
\end{equation*}
Matrix $L^{cp}_N$ is diagonalized by a discrete Fourier transform. It can be verified that its eigenvalues and eigenvectors are given by
\begin{align}
\lambda_0&=0 \label{eigen_cp1.eq} \\
\lambda_n&=N,\ \ \ \ \ n=1,2, \cdots, N-1 \label{eigen_cp2.eq} \\
u_{nm}&=\frac{1}{\sqrt{N}}e^{i2 \pi n m /N},\ \ \ \ \ n,m=0,1,\cdots,N-1\,. \label{eigen_cp3.eq}
\end{align}
\begin{proposition}
In a complete graph $\G=(\V,\E)$ with $N$ vertices, let $j,k\in \V$, $j\neq k$. The biharmonic distance between $j$ and $k$ is
\begin{align}
\label{cpBih:rst}
d_B(j,k)=\frac{\sqrt{2}}{N}\,.
\end{align}
\end{proposition}
\begin{IEEEproof}
By substituting the eigenvalues and eigenvectors in (\ref{eigen_cp1.eq}) - (\ref{eigen_cp3.eq}) into (\ref{biharmimg.eq}), we obtain
\begin{equation*}
\begin{split}
d^2_B(j,k)&=\sum^{N-1}_{n=1}\frac{|u_{nj}-u_{nk}|^2}{N^2}\\
&=\frac{1}{N^3}\sum^{N-1}_{n=1}4\sin^2\frac{(j-k)\pi n}{N}\\
&=\frac{2}{N^2}\,.
\end{split}
\end{equation*}
\end{IEEEproof}
Once we obtain the biharmonic distance between any vertices $j$ and $k$, we can derive the other related indices.
From~(\ref{cpBih:rst}), we derive the biharmonic Kirchhoff index for a complete graph with $N$ vertices.
\begin{align*}
D^2_B(\G)&=\frac{N(N-1)}{2}\cdot \frac{2}{N^2}=\frac{N-1}{N}\,.
\end{align*}
We also derive the biharmonic vertex index and biharmonic centrality for a complete graph,
\begin{align*}
D^2_B(j)&=(N-1)\cdot \frac{2}{N^2}=\frac{2(N-1)}{N^2}\,,\\
C_B(j)&=\frac{N^3}{2(N-1)}\,.
\end{align*}

Finally, we use the biharmonic distance and Theorems~\ref{hsojk_bih:them}, \ref{hsoj_bih:them}, and \ref{hsog_bih:them} to determine closed-form solutions for the three performance measures defined in Section~\ref{perfmeas:sec}.
\begin{theorem}
For a complete graph $\G$ with $N$ vertices, where the system dynamics are as given in (\ref{system:eqn}),
\begin{align*}
\Hso(j,k)&=\frac{1}{N^2}\,,\qquad j,k\in V,\,\, j\neq k\,;\\
\Hso(j)&=\frac{N-1}{2N^3}\,,\qquad j\in V\,;\\
\Hso(\G)&=\frac{N-1}{2N^2}\,.
\end{align*}
\end{theorem}


We recall that in a complete graph, the  total variance in a system with first-order noisy consensus dynamics (\ref{fosystem.eqn}) is $\Hfo(\G) \in O(1)$~\cite{YoScLe10}.
This is in contrast with $\Hso(\G)$ which is in $O(1/N)$.

\subsection{Star Graph}
We consider a star graph of order $N$, which consists of one hub and $N-1$ leaves. 
Its Laplacian matrix is
\begin{equation}
L^{star}_N=\left(
\begin{array}{ccccc}
  {N-1}&-1&{\cdots}&-1&-1  \\
  -1&1&{\cdots}&0&0 \\
  {\vdots}&{\vdots}&{\ddots}&{\vdots}&{\vdots} \\ 
-1&0&{\cdots}&1&0 \\   
  -1&0&{\cdots}&0&1 \\ 
\end{array}
\right) .
\end{equation}

Its eigenvalues and corresponding orthonormal eigenvectors are~\cite{Xu09},
\begin{align}
&\lambda_{0}=0\, \label{eigenva_star1.eq} \\
&\lambda_n=1\,,\quad n=1,2,\cdots,N-2\,,\label{eigenva_star2.eq} \\
&\lambda_{N-1}=N\,, \label{eigenva_star3.eq}
\end{align}
and
\begin{align}
u_{0}&=\frac{1}{\sqrt{N}}(1,1,1,\cdots,1,1,1)^\top\, \label{eigenvec_star1.eq} \\
u_{n}&=\frac{1}{\sqrt{n(n+1)}}(0,\underbrace{-1,\cdots,-1}_\text{$n$ },n,0,0,\cdots,0 )^\top, \nonumber \\
&\qquad n=1,2,\cdots,N-2\,,  \label{eigenvec_star2.eq}\\
u_{N-1}&=\frac{1}{\sqrt{N(N-1)}}(1-N,1,\cdots,1,1)^\top\,. \label{eigenvec_star3.eq}
\end{align}

We use these eigenvalues and eigenvectors to find the biharmonic distances between vertices in a star graph.
\begin{proposition} \label{starbiharm.prop}
In a star network $\G=(\V,\E)$ with vertex $0$ being the hub with degree $N-1$, and the remaining $N-1$ vertices as leaves, the biharmonic distance between the hub and a leaf is given by
\begin{equation}
\begin{split}
\label{stBih:res1}
d_B(0,j)&=\sqrt{\frac{N-1}{N}},\ \ \ \ \ j=1,2,\cdots,N-1,
\end{split}
\end{equation}
and the biharmonic distance between any two leaves is
\begin{equation}
\label{stBih:res2}
d_B(j,k)=\sqrt{2},\ \ \ \ \ \ j,k=1,2,\cdots,N-1; \quad j\neq k\,.
\end{equation}
\end{proposition}
\begin{IEEEproof}
The biharmonic distance between any two vertices $j,k\in \V$, $j\neq k$ is given by
\begin{equation}
\label{star_bihSum}
d^2_B(j,k)=\sum^{N-2}_{n=1}\frac{(u_{nj}-u_{nk})^2}{1^2}+\frac{(u_{N-1,j}-u_{N-1,k})^2}{N^2}.
\end{equation}
Substituting (\ref{eigenva_star1.eq}) - (\ref{eigenva_star3.eq}) and (\ref{eigenvec_star1.eq}) - (\ref{eigenvec_star3.eq}) into (\ref{star_bihSum}) yields the theorem.
\end{IEEEproof}

With these biharmonic distances, we easily obtain the biharmonic Kirchhoff index,
\begin{align*}
D^2_B(\G)&=(N-1)\frac{N-1}{N}+\frac{(N-1)(N-2)}{2}\cdot 2\nonumber\\
&=N^2-2N+\frac{1}{N}\,.
\end{align*}
The expressions for biharmonic vertex index and biharmonic centrality also immediately follow from the proposition, 
For the central vertex in a star graph,
\begin{align*}
D^2_B(0)&=(N-1)\cdot \frac{N-1}{N}=\frac{(N-1)^2}{N}\,,\\
C_B(0)&=\frac{N^2}{(N-1)^2}\,,
\end{align*}
and any leaf vertex $j$,
\begin{align*}
D^2_B(j)&=\frac{N-1}{N}+(N-2)\cdot 2=\frac{2N^2-3N-1}{N}\,,\\
C_B(j)&=\frac{N^2}{2N^2-3N-1}\,.
\end{align*}

Applying Proposition~\ref{starbiharm.prop} and Theorems~ \ref{hsojk_bih:them}, \ref{hsoj_bih:them}, and \ref{hsog_bih:them}, we obtain
closed-form solutions for the three steady-state variance performance measures.
\begin{theorem}
For a star graph $\G$ with $N$ vertices, where the system dynamics are as given in (\ref{system:eqn}), and where vertex $0$ is the hub,
\begin{align*}
\Hso(0,j)&=\frac{N-1}{2N}\,,\qquad j\neq 0\,;\\
\Hso(j,k)&=1\,,\qquad j\neq k;\,\, j,k\neq 0\,;\\
\Hso(0)&=\frac{N-1}{2N^3}\,;\\
\Hso(j)&=\frac{N^3-N^2-N-1}{2N^3}\,,\qquad j\neq 0\,;\\
\Hso(\G)&=\frac{N}{2}-1+\frac{1}{2N^2}\,.
\end{align*}
\end{theorem}
We recall that in an $N$-node star graph, the total variance for a system with first-order noisy consensus dynamics is $\Hfo(\G) \in O(N)$~\cite{YoScLe10},
and interestingly, in second order systems, the total variance is also in $O(N)$.

\subsection{Cycle}
The Laplacian of a cycle $C_N$ with $N$ vertices is given by 
\begin{equation*}
L^{cyc}_N=\left(
\begin{array}{ccccccc}
  2& -1&0&{\cdots}&0&0&-1  \\
  -1&2&-1&{\cdots}&0&0&0 \\ 
  {\vdots}&{\vdots}&{\vdots}&{\ddots}&{\vdots}&{\vdots}&{\vdots} \\ 
  0&0&0&{\cdots}&-1&2&-1 \\ 
  -1&0&0&{\cdots}&0&-1&2\\  
\end{array}
\right) \,.
\end{equation*}
$L^{cyc}_N$ is a circulant matrix. Therefore, its spectrum is given by a discrete Fourier transform. Let $\phi_n=\frac{n\pi}{N}$; the eigenvalues and eigenvectors of
$L^{cyc}_N$ are 
\begin{align}
\lambda_n&=2(1-\cos2\phi_n),n=0,1,2,\cdots,N-1 \label{per_eigen1.eq}\\
u_{nm}&=\frac{1}{\sqrt{N}}e^{i2m\phi_n}, n,m=0,1,\cdots,N-1. \label{per_eigen2.eq}
\end{align}

We use these eigenvalues and eigenvectors to determine the biharmonic distance.
\begin{proposition}
\label{cyc_bih_result}
In a cycle graph $\G = (\V,\E)$, let $j,k\in \V$, $k\leq j$ and $j-k=l$. Then, the biharmonic distance between $j$ and $k$ is
\begin{equation}
\label{cycBihRes:eqn}
d_B(j,k)=\sqrt{\frac{l^4}{12 N}-\frac{l^3}{6}+\frac{l^2 N}{12}-\frac{l^2}{6 N}+\frac{l}{6}}\,.
\end{equation}
\end{proposition}
The proof of the proposition is given in Appendix~\ref{ringProof:sec}.

Next, we calculate the derived indices  using biharmonic distances.
For a cycle $C_N$ with $N$ nodes, the biharmonic Kirchhoff index is
\begin{align*}
D^2_B(\G)&=\frac{1}{720}(N^5+10N^3-11N)\,.\nonumber
\end{align*}
For any vertex $j$ in a cycle, its biharmonic vertex index and biharmonic centrality are
\begin{align*}
D^2_B(j)&=\frac{1}{360}(N^4+10N^2-11)\,,\\
C_B(j)&=\frac{360N}{N^4+10N^2-11}\,.
\end{align*}

By applying Theorems~\ref{hsojk_bih:them}, \ref{hsoj_bih:them}, and \ref{hsog_bih:them}, along with Proposition~\ref{cyc_bih_result},
we obtain closed-form solutions for the steady-state variance performance measures.
\begin{theorem}
For a cycle graph $\G$ with $N$ vertices where the dynamics are given by  (\ref{system:eqn}),
\begin{align}
\Hso(j,k)&=\frac{l^4}{24 N}-\frac{l^3}{12}+\frac{l^2 N}{24}-\frac{l^2}{12 N}+\frac{l}{12}\,,\nonumber\\
&\qquad\text{For } j,k\in \V,\quad k\leq j \text{ and } j-k=l\,;\\
\Hso(j)&=\frac{1}{1440}\left(N^3+10N-\frac{11}{N}\right)\,,\quad j\in \V\,;\\
\Hso(\G)&=\frac{1}{1440}\left(N^4+10N^2-11\right)\,.
\end{align}
\end{theorem}

To give some examples for $\Hso(j,k)$ in a cycle of $N$ vertices, it holds that $\Hso(0,1)=\frac{1}{24}(N-1/N)$. For a even $N$, $\Hso(0,N/2)=\frac{1}{384}N(N^2+8)$.

To compare with the first-order consensus dynamics, we recall that in a cycle graph with $N$ vertices, 
$\Hfo(\G) \in O(N^2)$~\cite{BaJoMiPa12}, whereas in second-order systems $\Hso(\G) \in O(N^4)$.

\subsection{Path}
We consider a path graph $P_N$ with $N$ vertices.  Let the vertices be numbered $0, 1, \ldots, N-1$.
The Laplacian matrix of $P_N$ assumes the form
\begin{equation*}
L^{path}_N=\left(
\begin{array}{ccccccc}
  1& -1&0&{\cdots}&0&0&0  \\
  -1&2&-1&{\cdots}&0&0&0 \\ 
  {\vdots}&{\vdots}&{\vdots}&{\ddots}&{\vdots}&{\vdots}&{\vdots} \\ 
  0&0&0&{\cdots}&-1&2&-1 \\ 
  0&0&0&{\cdots}&0&-1&1\\  
\end{array}
\right) .
\end{equation*}
The eigenvalues and eigenvectors of $L_N^{path}$ are~\cite{TzWu00}.
\begin{align}
\lambda_n=&2(1-\cos\phi_n),\ \ \ \ \ n=0,1,2,\cdots,N-1 \label{free_eigen1.eq} \\
u_{0m}=&\frac{1}{\sqrt{N}}, \ \ \ \ \ m=0,1,\cdots N-1 \label{free_eigen2.eq}\\
u_{nm}=&\sqrt{\frac{2}{N}}\cos(m+1/2\phi_n), \ \ \ \ \nonumber\\ &n=1,2,\cdots,N-1,m=0,1,\cdots N-1 \label{free_eigen3.eq}
\end{align}
where $\phi_n=n\pi/N$.

We use (\ref{free_eigen1.eq}) - (\ref{free_eigen3.eq}) to determine the biharmonic distance between two vertices in a path.
\begin{proposition}
\label{path_bih_result}
In a path graph $\G= (\V,\E)$ with $N$ vertices, 
the biharmonic distance between two vertices ${j, k\in  \{0,1, \ldots, N-1\}}$, $k<j$, is 
\begin{small}
\begin{equation}
\label{pathBihResult:eqn}
\begin{split}
d_B(j,k)&=\Bigg(\frac{j}{6}
+\frac{j^2}{2}
-\frac{j^2}{4N}+\frac{j^3}
{3}-\frac{j^3}{2N}-\frac{j^4}{4N}\\
&-\frac{k}{6}-jk+\frac{jk}{2N}+\frac{j^2k}{2N}+\frac{k^2}{2}-\frac{k^2}{4N}\\&-jk^2+\frac{jk^2}{2N}+\frac{j^2 k^2}{2N}+\frac{2k^3}{3}-\frac{k^3}{2N}-\frac{k^4}{4N}\Bigg)^{\frac{1}{2}}\,.
\end{split}
\end{equation}
\end{small}
\end{proposition}
The proof of Proposition~\ref{path_bih_result} is given in Appendix \ref{pathProof:sec}.

We next use Proposition~\ref{path_bih_result} to derive the biharmonic Kirchhoff index for a path with $N$ nodes,
\begin{align*}
D^2_B(\G)&=\frac{1}{180}(2N^5+5N^3-7N)\,.\nonumber
\end{align*}
We can also derive the biharmonic vertex index and biharmonic centrality for a node $j$,
\begin{align*}
 D^2_B(j)=&\frac{1}{30}(N^4-10j(j+1)N^2\nonumber\\
 &+10j(2j+1)(j+1)N-10j^2(j+1)^2-1)\,,\\
  C_B(j)=&\textstyle\frac{30N}{N^4-10j(j+1)N^2+10j(2j+1)(j+1)N-10j^2(j+1)^2-1}\,.
\end{align*}
%
%
Finally, we present the following theorem that gives the steady-state variance performance measures for $P_N$.
This theorem follows directly from Proposition~\ref{path_bih_result} and Theorems~\ref{hsojk_bih:them}, \ref{hsoj_bih:them}, and \ref{hsog_bih:them}.
\begin{theorem}
Let  $\G = (\V,\E)$ be a path graph with ${\V = \{0,1, \ldots, N-1\}}$ and with the dynamics (\ref{system:eqn}). Let $j, k \in \V$ with $k < j$.
Then,
\begin{align}
&\Hso(j,k)=\frac{1}{2}\Bigg(\frac{j}{6}
+\frac{j^2}{2}
-\frac{j^2}{4N}+\frac{j^3}
{3}-\frac{j^3}{2N}-\frac{j^4}{4N}\nonumber\\
&~~~-\frac{k}{6}-jk+\frac{jk}{2N}+\frac{j^2k}{2N}+\frac{k^2}{2}-\frac{k^2}{4N}\nonumber\\
&~~~-jk^2+\frac{jk^2}{2N}+\frac{j^2 k^2}{2N}+\frac{2k^3}{3}-\frac{k^3}{2N}-\frac{k^4}{4N}\Bigg)\,,\\
&\Hso(j)=\frac{1}{360N}(4N^4-(60j^2+60j+5)N^2\nonumber\\
 &~~~+60j(2j+1)(j+1)N-60j^2(j+1)^2+1)\\
&\Hso(\G)=\frac{1}{360}\left(2N^4+5N^2-7\right)\,.
\end{align}
\end{theorem}

To give some examples for $\Hso(j,k)$ and $\Hso(j)$, we note that $\Hso(0,N-1)=\frac{1}{24}N(N^2-1)$ and $\Hso(0)=\frac{1}{90}N(N^2-5/4)+\frac{1}{360N}$. For a even $N$, $\Hso(N/2,N-1)=\frac{5}{384}N(N-2/5)(N-2)$ and $\Hso(N/2)=\frac{1}{1440}N(N^2+40)+\frac{1}{360N}$. 

We recall that in a $N$-vertex path graph with first order noisy consensus dynamics, the total variance is $\Hfo = O(N^2)$~\cite{YoScLe10}.
This is in contrast with the second order system, which has total variance in $O(N^4)$.

\section{Numerical Examples}
\label{num:sec}

In this section, we give  numerical examples of the biharmonic and resistance distances in several graphs.

\begin{figure}[htbp]
\begin{center}
\subfigure{\includegraphics[width=0.48\linewidth]{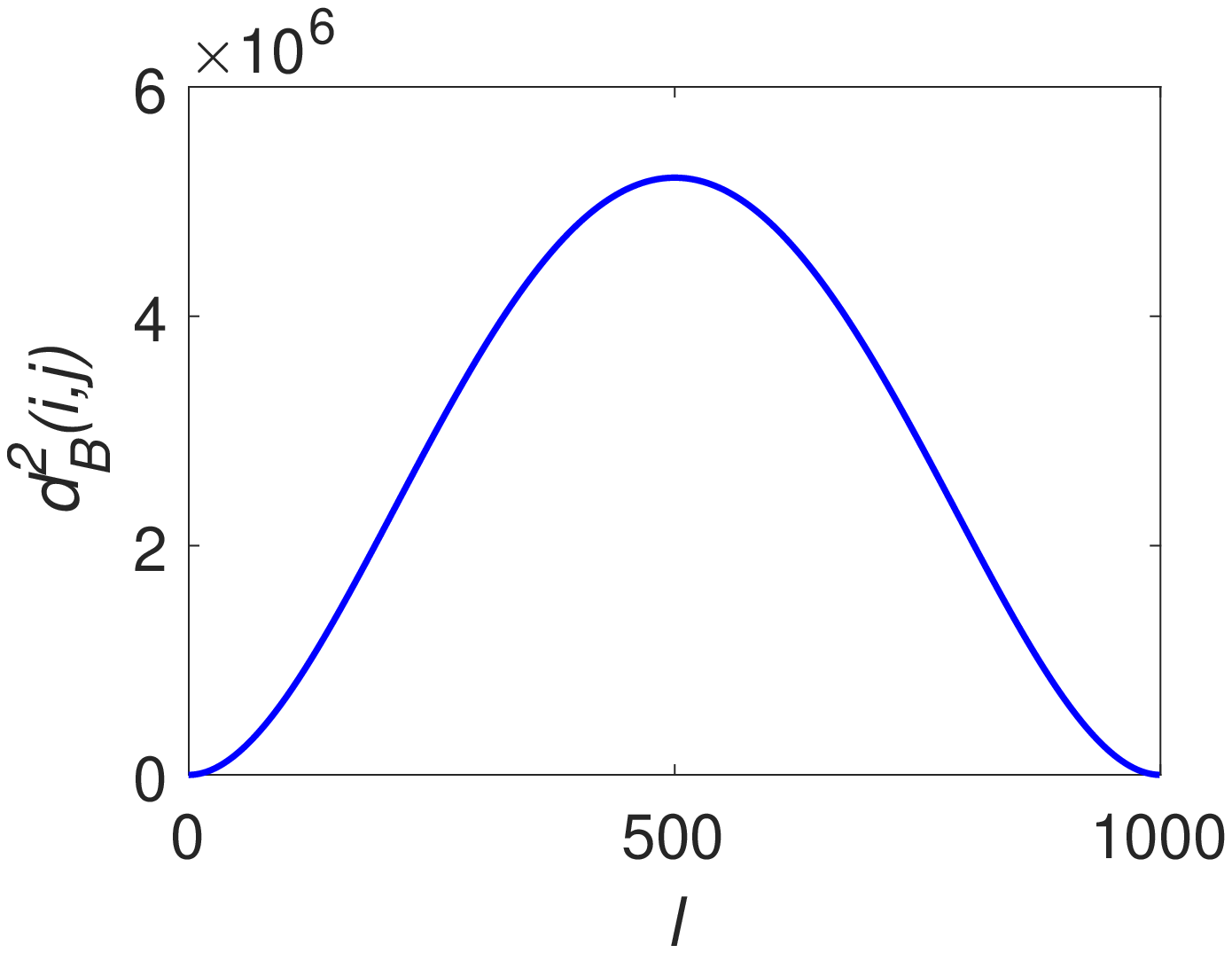}}
\subfigure{\includegraphics[width=0.48\linewidth]{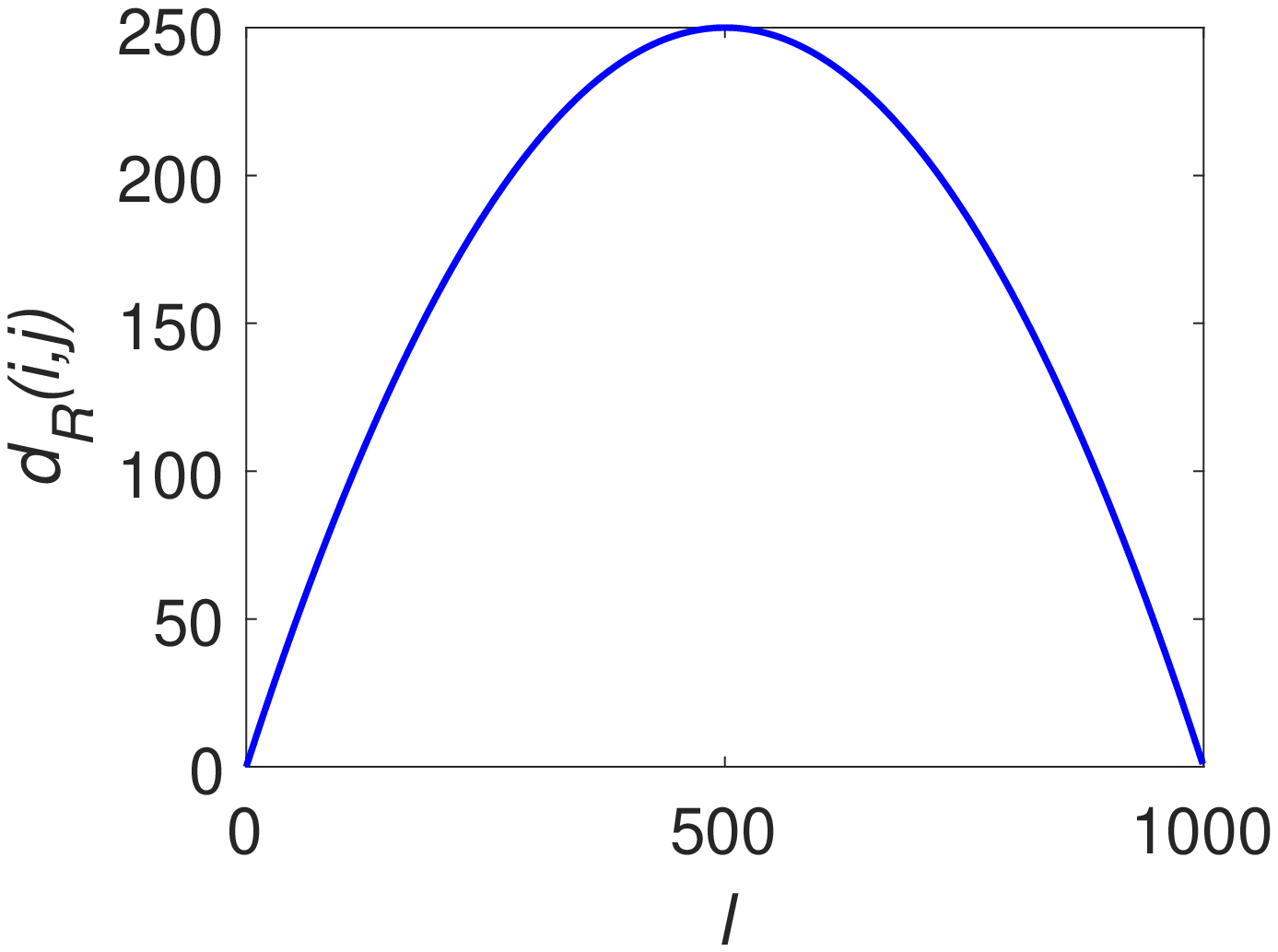}}
\caption{The squared biharmonic distance $d^2_B(j,k)$ and resistance distance $d_R(j,k)$ between two vertices $j,k$ with $l=j-k$ in a cycle of $1000$ vertices.}
\label{cycDistCpr.fig}
\end{center}
\end{figure}

Figure \ref{cycDistCpr.fig} shows the square of biharmonic distance and the resistance distance in a cycle of $1000$ vertices. Specifically we plot both the distances between vertices $j$ and $k$ where $k\leq j$, as a function of $l = j-k$. The biharmonic distances are obtained using (\ref{cycBihRes:eqn}).
The figure shows that the square of biharmonic distance and the resistance distance grow at different rates in a cycle, as a function of graph distance, while the vertices that have the largest graph distance have both the largest squared biharmonic distance and resistance distance.

\begin{figure}[htbp]
\begin{center}
\subfigure[$k=0$]{\includegraphics[width=0.48\linewidth]{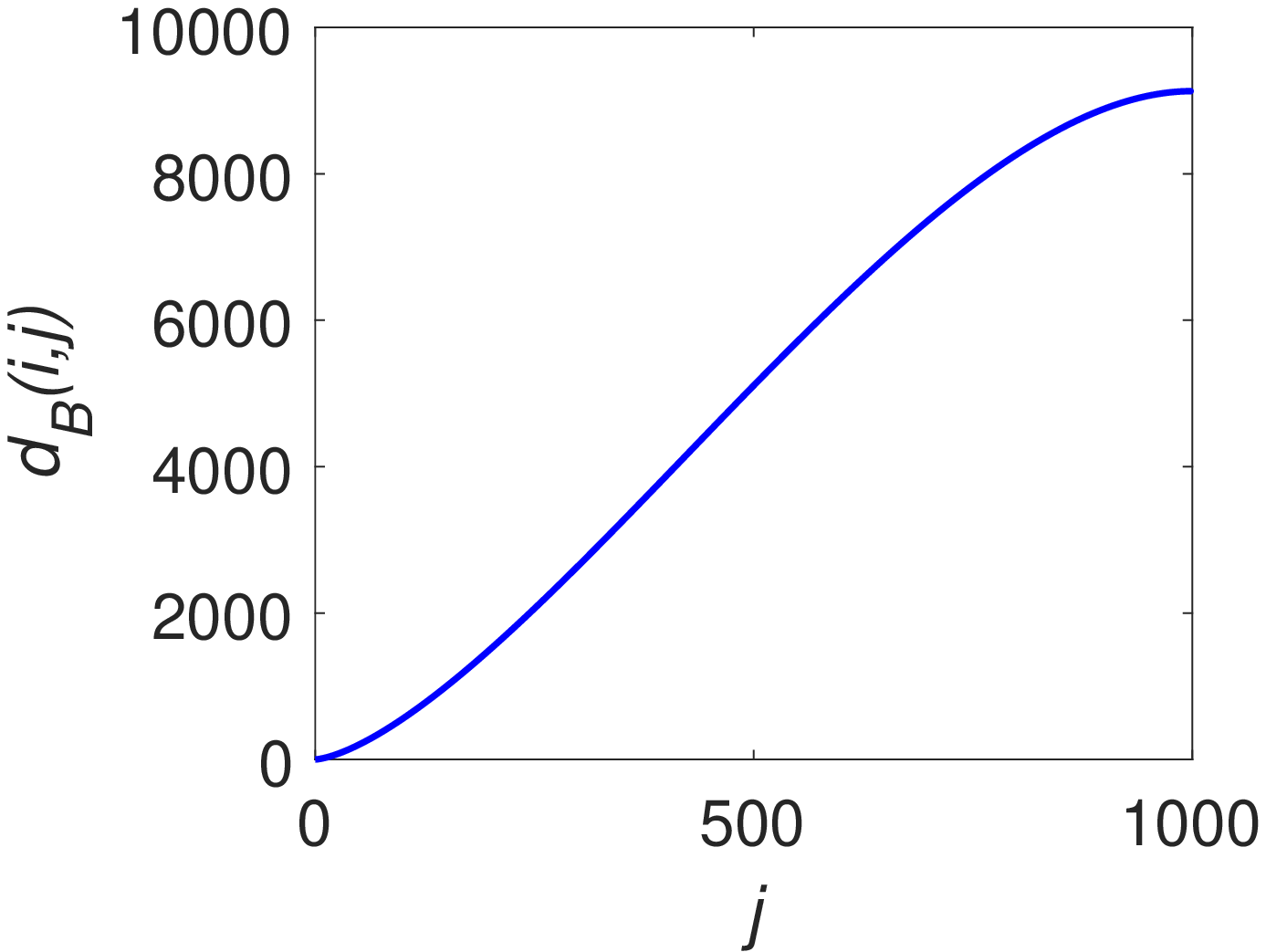}}
\subfigure[$k=500$]{\includegraphics[width=0.48\linewidth]{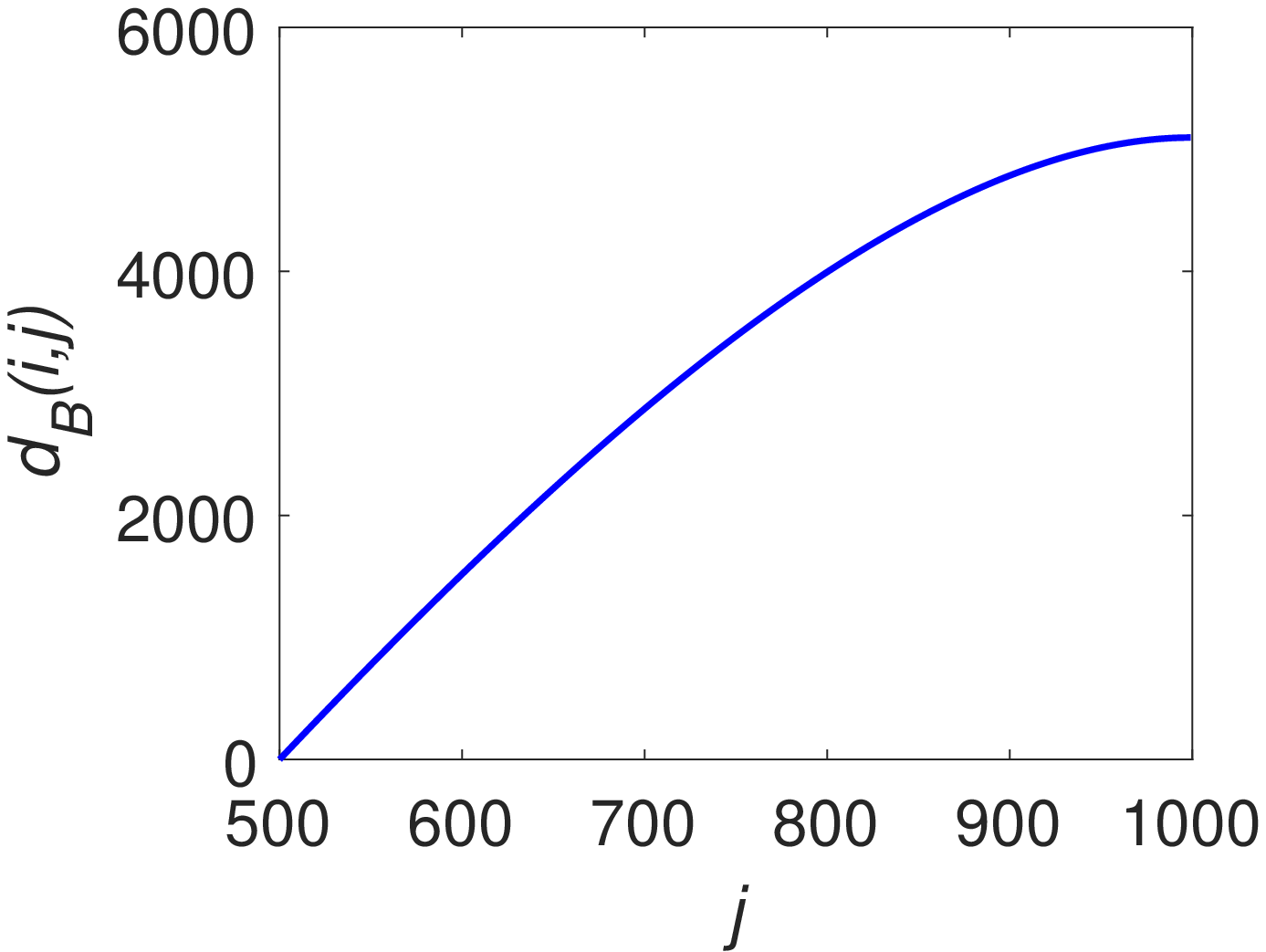}}
\caption{Biharmonic distance $d_B(j,k)$ between two vertices $j,k$ with $l=j-k$ in a path of $1000$ vertices.} 
\label{pathBihDist:fig}
\end{center}
\end{figure}

Figure \ref{pathBihDist:fig} gives the biharmonic distances in a path graph. In particular, we show the biharmonic distances between vertices $j$ and $k$ where $k\leq j$. We only show two  cases, $k=0$ and $k=500$. The biharmonic distances are calculated using (\ref{pathBihResult:eqn}). For a given $k$, $d_B(j,k)$ grows slower near the ends of the path and faster around the middle of the path. In addition, since for even $N$, $d_B(0,N/2-1)=d_B(N/2,N-1)$; we observe that $d_B(0,N/2-1)+d_B(N/2-1,N/2)+d_B(N/2,N-1)>d_B(0,N-1)$ in this example. This is in contrast with resistance distance (and identically graph distance), where $d_R(0,N/2-1)+d_R(N/2-1,N/2)+d_R(N/2,N-1)=d_R(0,N-1)$.

\begin{figure}[htbp]
\begin{center}
\subfigure{\includegraphics[width=0.48\linewidth]{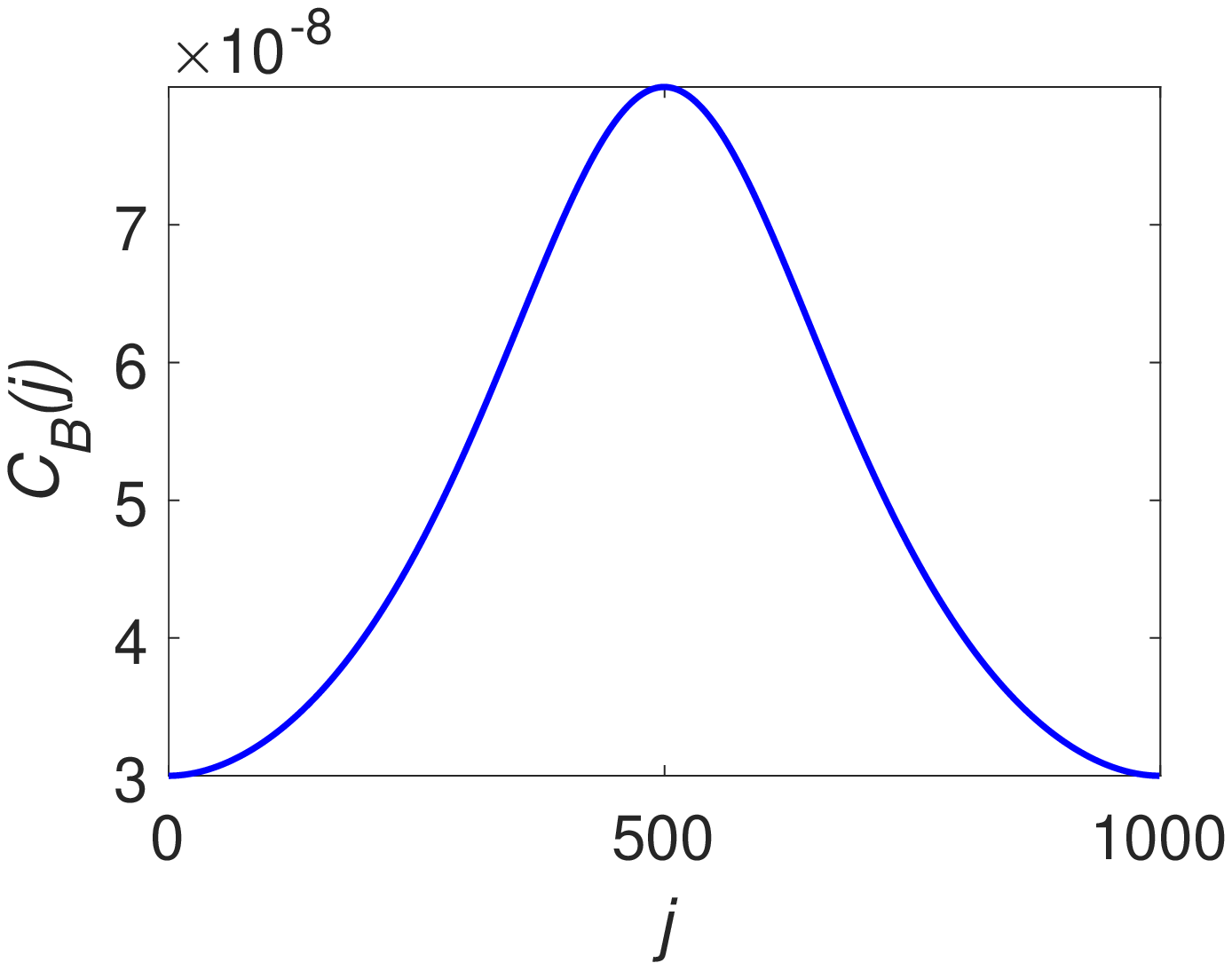}}
\subfigure{\includegraphics[width=0.48\linewidth]{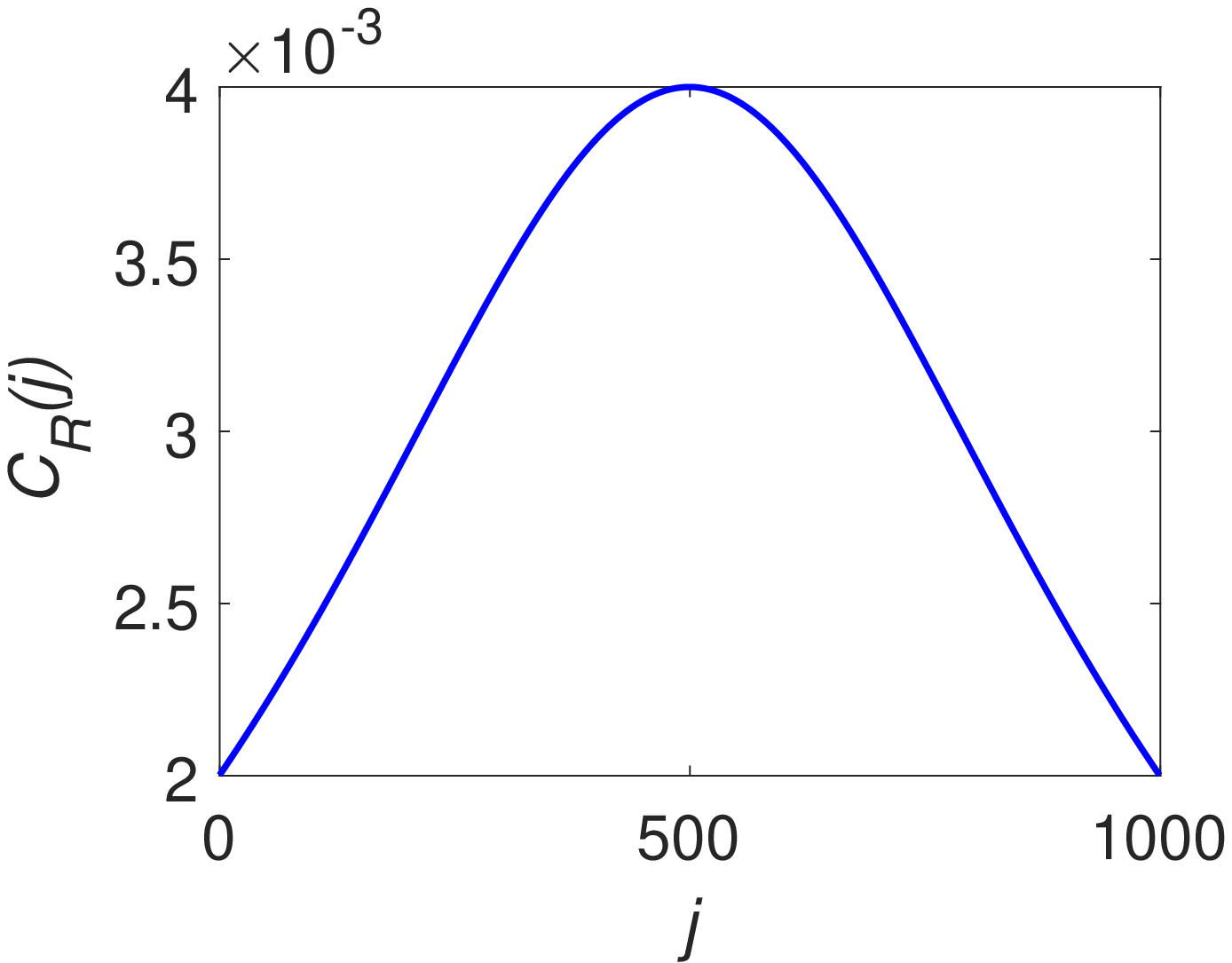}}
\caption{Biharmonic centrality and information centrality in a path of $1000$ vertices.}
\label{CentjLine.fig}
\end{center}
\end{figure}

Figure \ref{CentjLine.fig} compares biharmonic centrality and information centrality in a path with 1000 vertices.
Both curves are bell-like and the node in the middle has the largest centrality. The difference is that biharmonic distance  distinguishes the center nodes better, as illustrated by the figure.

\begin{figure} 
\centering
\includegraphics[width=0.9\linewidth, height=0.5\linewidth]{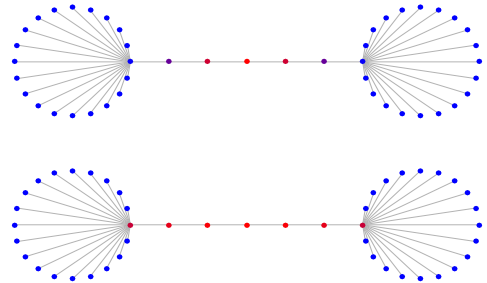}
\caption{Biharmonic centrality and information centrality in a starry-line graph.} \label{starryLine.fig}
\end{figure}
The next example is a starry-line graph, composed of  two $20$-vertex star graphs connected by a path of $5$ vertices. Figure~\ref{starryLine.fig} shows the biharmonic centrality (above) and information centrality (below) in the graph. Vertices are colored according to their centrality in the network. Red vertices have the largest centrality and blue vertices have smallest centrality. The figure shows that the biharmonic centrality distinguishes the center of the line from other vertices on the line, while these vertices have comparable information centralities.

\begin{figure}[htbp]
\begin{center}
\subfigure[Biharmonic embedding and biharmonic centrality.]{\label{bihCent:fig}\includegraphics[width=.7 \linewidth]{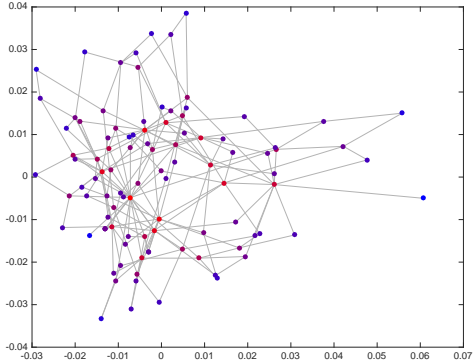}}\\
\subfigure[Resistance embedding and information centrality.]{\label{resCent:fig}\includegraphics[width=.7 \linewidth]{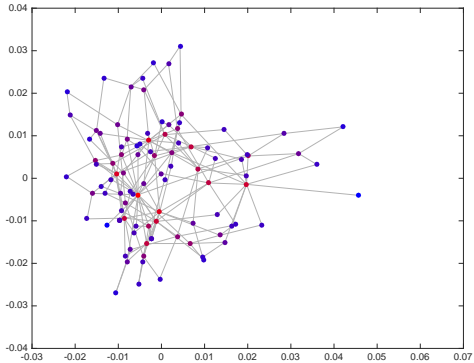}}
\caption{Embeddings and centralities of a $100$-vertex BA network}
\label{cent:fig}
\end{center}
\end{figure}

Figure \ref{cent:fig} shows the first two principle components of the biharmonic embedding as well as the biharmonic embedding of a Barab\'{a}si–-Albert network with $100$ nodes. 
We observe that the biharmonic embedding stretches the edges out a bit more than the resistance embedding. In fact, by reviewing their definitions, we  observe that the normalized components in PCA for these two embeddings are the same; the differences are the variances of the components.


\section{Conclusion} \label{conclusion:sec}
We have investigated the performance of undirected networks with second-order consensus dynamics with stochastic disturbances. We have established the connection between second-order network performance measures and the biharmornic distances in the communication graph. We introduced the notions of a Kirchhoff index and vertex centrality based on biharmonic distance to further help us describe the behavior of second-order consensus dynamics, and we derived closed-form expressions for the performance measures for complete graphs, star graphs, cycles, and paths. 
Future work should include the study of additional properties of biharmonic distances,
as well as analysis of the steady-state variance performance measures   in more general networks, including random networks and real-world networks. 

\appendices
\section{Trigonometric Identities}
\label{eqns:sec}
We use the notation  $\phi_n=\frac{n\pi}{N}$. We next introduce the following identities.
\begin{small}
\begin{align}
\label{idty1}
G_N(1)&=\frac{1}{N}\sum^{N-1}_{n=1}\frac{1-\cos2\phi_n}{(1-\cos2\phi_n)^2}\nonumber\\
&=\frac{1}{2N}\sum^{N-1}_{n=1}\frac{1}{\sin^2 \phi_n}=\frac{N}{6}-\frac{1}{6 N}
\end{align}
\begin{align}
\label{idty2}
G_N(2)&=\frac{1}{N}\sum^{N-1}_{n=1}\frac{1-\cos4\phi_n}{(1-\cos2\phi_n)^2}\nonumber\\
&=\frac{1}{2N}\sum^{N-1}_{n=1}\frac{\sin^2 2\phi_n}{\sin^4\phi_n}=\frac{2}{N}\sum^{N-1}_{n=1}\frac{\cos^2\phi_n}{\sin^2\phi_n}\nonumber\\
&=\frac{2N}{3}-2+\frac{4}{3N}
\end{align}
\begin{align}
\label{idty3}
F_N(1)&=\frac{1}{N}\sum^{N-1}_{n=1}\frac{1-\cos\phi_n}{(1-\cos\phi_n)^2}\nonumber\\
&=\frac{1}{2N}\sum^{N-1}_{n=1}\frac{1}{\sin^2 \phi_n/2}=\frac{N}{3}-\frac{1}{3N}
\end{align}
\begin{align}
\label{idty4}
F_N(2)&=\frac{1}{N}\sum^{N-1}_{n=1}\frac{1-\cos2\phi_n}{(1-\cos\phi_n)^2}\nonumber\\
&=\frac{1}{2N}\sum^{N-1}_{n=1}\frac{\sin^2\phi_n}{\sin^4\phi_n/2}=\frac{2}{N}\sum^{N-1}_{n=1}\frac{\cos^2\phi_n/2}{\sin^2\phi_n/2}\nonumber\\
&=\frac{4N}{3}-2+\frac{2}{3N}
\end{align}
\end{small}

\section{Proof of Proposition \ref{cyc_bih_result}}
\label{ringProof:sec}
\begin{IEEEproof}
We note that $i$ denotes the imaginary unit in this proof.

Substituting (\ref{per_eigen1.eq}) and (\ref{per_eigen2.eq}) into Definition \ref{biharmonic:def}, we obtain
\begin{small}
\begin{align}
\label{dbound}
d^2_B(j,k)&=\frac{1}{N}\sum_{n=1}^{N-1}\frac{|e^{i2j\phi_n}-e^{i2k\phi_n}|^2}{4(1-\cos2\phi_n)^2}=\frac{1}{2}G_N(j-k),
\end{align}
\end{small}
where
\begin{equation*}
G_N(l)=\frac{1}{N}\sum_{n=1}^{N-1}\frac{1-\cos(2l\phi_n)}{(1-\cos2\phi_n)^2}\,.
\end{equation*}
Without loss of generality, we assume $0 \leqslant l\leqslant 2N$.

In order to simplify $G_N(l)$, we give two equivalent expressions for the real part of the following sum
\begin{align}
\label{hnl:def}
H_N(l)=\frac{1}{N}\sum_{n=1}^{N-1}\frac{1-e^{2i l \phi_n}}{(1-e^{2i \phi_n})^2}\,.
\end{align}
The first expression is 
\begin{small}
\begin{align}\label{recursion}
&\mathrm{Re}\big(H_N(l)\big)\nonumber\\
&=\frac{1}{4N} \sum_{n=1}^{N-1}\bigg(\frac{1-\cos 2 l\phi_n}{(1-\cos\phi_n)^2}-\frac{2-2\cos2(l-1)\phi_n}{{(1-\cos\phi_n)^2}}\nonumber\\
&\quad +\frac{1-\cos2(l-2)\phi_n}{{(1-\cos\phi_n)^2}}-\frac{1-\cos4\phi_n}{{(1-\cos\phi_n)^2}}+\frac{2-2\cos 2\phi_n}{(1-\cos\phi_n)^2} \bigg)\nonumber\\
&=\frac{1}{4}\bigg(G_N(l)-2G_N(l-1)\nonumber\\
&\qquad\qquad +G_N(l-2)-G_N(2)+2G_N(1)\bigg)\,.
\end{align}
\end{small}
We note that $G_N(0)=0$. Let $K_N(l)=G_N(l)-G_N(l-1)$. We rewrite (\ref{recursion}) for the sake of conciseness in future derivation as
\begin{small}
\begin{align}\label{r1}
&\mathrm{Re}( H_N(l))\nonumber\\
&=\frac{1}{4}\Big(\left[\left(G_N(l)-G_N(l-1)\big)-\big(G_N(l-1)-G_N(l-2)\right)\right]\nonumber\\
&\ \ \ \ \ -\left[\left(G_N(2)-G_N(1)\right)-\left(G_N(1)-G_N(0)\right)\right]\Big)\nonumber\\
&=\frac{1}{4}\left[\left(K_N(l)-K_N(l-1)\right)-\left(K_N(2)-K_N(1)\right)\right]
\end{align}
\end{small}
Next, we use the summation formula $\sum^{n-1}_{j=0}x^j=\frac{1-x^n}{1-x}$ to expand (\ref{hnl:def})
\begin{small}
\begin{align}
H_N(l)=&\frac{1}{N} \sum_{n=1}^{N-1}\frac{1-e^{2il\phi_n}}{1-e^{2i\phi_n}}\frac{1}{1-e^{2i\phi_n}}\nonumber\\
=&\frac{1}{N} 
\sum_{n=1}^{N-1}
\sum_{l'=0}^{l-1}
\bigg(\frac{e^{2il'\phi_n}}{1-e^{2i\phi_n}}-\frac{1}{1-e^{2i\phi_n}}+\frac{1}{1-e^{2i\phi_n}}\bigg)\nonumber\\
=&\frac{1}{N} 
\sum_{n=1}^{N-1}
\bigg(\sum_{l'=2}^{l-1}
\frac{e^{2il'\phi_n}-1}{1-e^{2i\phi_n}}-1\bigg)+\frac{1}{N} 
\sum_{n=1}^{N-1}
\sum_{l'=0}^{l-1}
\frac{1}{1-e^{2i\phi_n}}\nonumber\\
=&-\frac{1}{N} 
\sum_{n=1}^{N-1}
\sum_{l'=2}^{l-1}
\sum_{l''=1}^{l'-1}
e^{2il''\phi_n}-\frac{1}{N}\sum^{N-1}_{n=1}1\nonumber\\
&-\frac{1}{N}\sum_{n=1}^{N-1}
\sum_{l'=2}^{l-1}1+\frac{1}{N} 
\sum_{n=1}^{N-1}
\sum_{l'=0}^{l-1}
\frac{1}{1-e^{2i\phi_n}}\,. \label{HN.eq}
\end{align}
\end{small}
The triple summation in the last equality can be simplified by carrying out the summation over $n$ first, 
\begin{align*}
E_1&\equiv-\frac{1}{N} 
\sum_{n=1}^{N-1}
\sum_{l'=2}^{l-1}
\sum_{l''=1}^{l'-1}
e^{2il''\phi_n}\nonumber\\
&=-\frac{1}{N}
\sum_{l'=2}^{l-1}\sum_{l''=1}^{l'-1}\left(\frac{1-e^{i2\pi l''}}{1-e^{i \pi y'' /N}}-1\right)\nonumber\\
&=\frac{1}{N}
\sum_{l'=2}^{l-1}\sum_{l''=1}^{l'-1}1,
\end{align*}
where last equality is obtained by applying $e^{i2\pi l''}=1$ for $l''\in {\mathbb{Z}}$.

Using the fact that ${\mathrm{Re}}\left(1/(1-e^{i\theta})\right)=1/2$, $0<\theta<2 \pi$, the real part of the fourth term in (\ref{HN.eq}) is
\begin{equation*}
\mathrm{Re} (E_4)=\mathrm{Re}\Bigg(\frac{1}{N} 
\sum_{n=1}^{N-1}
\sum_{l'=0}^{l-1}
\frac{1}{1-e^{i\phi_n}}\Bigg)=\frac{(N-1)l}{2N}\,.
\end{equation*}
Therefore,
 \begin{equation}\label{r2}
\mathrm{Re}(H_N(l))=\frac{l^2}{2  N}-\frac{l}{N}-\frac{l}{2}+1\,.
\end{equation}

Let $X_N(l)=4\mathrm{Re} (H_N(l))$. From the equivalence of (\ref{r1}) and (\ref{r2}), we derive
\begin{small}
\begin{equation*}
\begin{split}
X_N(l)= \Big(K_N(l)-K_N(l-1)\Big)-\Big(K_N(2)-K_N(1)\Big).
\end{split}
\end{equation*}
\end{small}
This recursive equation can be solved to give
\begin{align*}
K_N(l)&=G_N(l)-G_N(l-1)\\
&=Y_N(l) +(l-1)G_N(2)-(2l-3)G_N(1),
\end{align*}
and
\begin{equation}
\label{gnl:eqn}
G_N(l)=Z_N(l) +\bigg(\frac{l^2}{2}-\frac{l}{2}\bigg)G_N(2)-(l^2-2l)G_N(1),
\end{equation}
where
\begin{equation*}
\begin{split}
Y_N(l)=& \sum^{l}_{j=2}X_N(j) \quad\text{ and }\\
Z_N(l)=&\sum^{l}_{j=2}Y_N(j)
=\sum^{l}_{j=2}\sum^{j}_{k=2}X_N(k)\,.
\end{split}
\end{equation*}

Substituting (\ref{idty1}), (\ref{idty2}), and (\ref{r2}) into (\ref{gnl:eqn}), we finally obtain the result for $G_N(l)$ as
\begin{equation*}
G_N(l)=\frac{l^4}{6 N}-\frac{l^3}{3}+\frac{l^2 N}{6}-\frac{l^2}{3 N}+\frac{l}{3}\,.
\end{equation*}
Plugging this value into (\ref{dbound}) generates the result in Proposition \ref{cyc_bih_result}.
\end{IEEEproof}

\section{Proof of Proposition ~\ref{path_bih_result}}
\label{pathProof:sec}
\begin{IEEEproof}
We note that $i$ denotes the imaginary unit in this proof.

By definition,  the biharmonic distance between $j$ and $k$, $j\leq k$ is
\begin{small}
\begin{align}\label{dfree}
d^2_B(j,k)=&\frac{1}{N}\sum_{n=1}^{N-1}\frac{[\cos(j+\frac{1}{2})\phi_n-\cos(k+\frac{1}{2})\phi_n]^2}{2(1-\cos\phi_n)^2}\nonumber\\
=&\frac{1}{2}\Big(F_N(j+k+1)+F_N(j-k)\nonumber\\&-\frac{1}{2}F_N(2j+1)-\frac{1}{2}F_N(2k+1)\Big)
\end{align}
where 
\begin{equation*}
F_N(l)=\frac{1}{N}\sum_{n=1}^{N-1}\frac{1-\cos l \phi_n}{(1-\cos\phi_n)^2}\,.
\end{equation*}
\end{small}
Next, we calculate the real part of the following sum in two different ways
\begin{equation}
T_N(l)= \frac{1}{N}\sum_{n=1}^{N-1} \frac{1-e^{il\phi_n}}{(1-e^{i\phi})^2}.
\end{equation}

First, let $E_N(l) = F_N(l)-F_N(l-1)$.  We obtain
\begin{small}
\begin{align}\label{retn1}
&\mathrm{Re}\big(T_N(l)\big)
=\frac{1}{4}\bigg(F_N(l)-2F_N(l-1)+F_N(l-2)\nonumber\\
&\qquad\qquad -F_N(2)+2F_N(1)\bigg)\nonumber\\
&=\frac{1}{4}\Big[\big(E_N(l)-E_N(l-1)\big)-\big(E_N(2)-E_N(1)\big)\Big]
\end{align}
\end{small}
Second, we use the summation formula $\sum^{n-1}_{j=0}x^j=\frac{1-x^n}{1-x}$ and derive
\begin{small}
\begin{align}
T_N(l)&=\frac{1}{N} \sum_{n=1}^{N-1}\frac{1-e^{il\phi_n}}{1-e^{i\phi_n}}\frac{1}{1-e^{i\phi_n}}\nonumber\\
=&\frac{1}{N} 
\sum_{n=1}^{N-1}
\sum_{l'=0}^{l-1}
\bigg(\frac{e^{il'\phi_n}}{1-e^{i\phi_n}}-\frac{1}{1-e^{i\phi_n}}+\frac{1}{1-e^{i\phi_n}}\bigg)\nonumber\\
=&\frac{1}{N} 
\sum_{n=1}^{N-1}
\bigg(\sum_{l'=2}^{l-1}
\frac{e^{il'\phi_n}-1}{1-e^{i\phi_n}}-1\bigg)+\frac{1}{N} 
\sum_{n=1}^{N-1}
\sum_{l'=0}^{l-1}
\frac{1}{1-e^{i\phi_n}}\nonumber\\
=&-\frac{1}{N} 
\sum_{n=1}^{N-1}
\sum_{l'=2}^{l-1}
\sum_{l''=1}^{l'-1}
e^{il''\phi_n}-\frac{1}{N}\sum^{N-1}_{n=1}1\nonumber\\
&-\frac{1}{N}\sum_{n=1}^{N-1}
\sum_{l'=2}^{l-1}1+\frac{1}{N} 
\sum_{n=1}^{N-1}
\sum_{l'=0}^{l-1}
\frac{1}{1-e^{i\phi_n}}\,. \label{TN.eq}
\end{align}
\end{small}
Again, we change the order of summation over $n$, $l'$ and $l''$ to simplify the first term in (\ref{TN.eq}),
\begin{small}
\begin{equation*}
E'_1\equiv-\frac{1}{N} 
\sum_{n=1}^{N-1}
\sum_{l'=2}^{l-1}
\sum_{l''=1}^{l'-1}
e^{il''\phi_n}=-\frac{1}{N}
\sum_{l'=2}^{l-1}\sum_{l''=1}^{l'-1}\bigg[\frac{1-(-1)^{l''}}{1-e^{i \pi y'' /N}}-1\bigg].
\end{equation*}
\end{small}
The real part of $E'_1$ and the  fourth term in (\ref{TN.eq}), denoted $E'_4$, are
\begin{equation*}
\mathrm{Re} (E'_1)=\frac{1}{8N}\left[2l^2-8l+7+(-1)^l\right],
\end{equation*}
\begin{equation*}
\mathrm{Re} (E'_4)=\mathrm{Re}\Bigg(\frac{1}{N} 
\sum_{n=1}^{N-1}
\sum_{l'=0}^{l-1}
\frac{1}{1-e^{i\phi_n}}\Bigg)=\frac{(N-1)l}{2N}.
\end{equation*}
Hence,
\begin{equation}\label{retn2}
\mathrm{Re} (T_N(l))=\frac{-4 N (l-2)+2 l^2-4l+(-1)^l-1}{8 N}.
\end{equation}
Equating (\ref{retn1}) and (\ref{retn2}) leads to
\begin{small}
\begin{align}
\label{fSumForm}
F_N(l)=&\sum^{l}_{i=2}\sum^{i}_{j=2}4\mathrm{Re} (T_N(l)) +\bigg(\frac{l^2}{2}-\frac{l}{2}\bigg)F_N(2)\nonumber\\
&-(l^2-2l)F_N(1)\,.
\end{align}
\end{small}
Next, we evaluate $F_N(1)$ and $F_N(2)$
\begin{small}
\begin{align}
F_N(1)&=\frac{1}{N}\sum^{N-1}_{n=1}\frac{1-\cos\phi_n}{(1-\cos\phi_n)^2}=\frac{1}{2N}\sum^{N-1}_{n=1}\frac{1}{\sin^2 \phi_n/2}\,,\\
F_N(2)&=\frac{1}{N}\sum^{N-1}_{n=1}\frac{1-\cos2\phi_n}{(1-\cos\phi_n)^2}=\frac{2}{N}\sum^{N-1}_{n=1}\frac{\cos^2\phi_n/2}{\sin^2\phi_n/2}\,.
\end{align}
\end{small}

For $F_N(1)$, we start by expanding the expression  $
\sum^{2N-1}_{n=1} 1/(\sin^2\frac{n\pi}{2N})$. Since $\sum^{N-1}_{n=1}1/(\sin^2\frac{n\pi}{N})=\frac{ N^2}{3}-\frac{1}{3}$, we derive
\begin{small}
\begin{align*}
&\sum^{2N-1}_{n=1}\frac{1}{\sin^2\frac{n\pi}{2N}}=\frac{ 4N^2}{3}-\frac{1}{3}\\
&=\frac{1}{\sin^2\frac{\pi}{2N}}+\frac{1}{\sin^2\frac{2\pi}{2N}}+\cdots+\frac{1}{\sin^2\frac{(N-1)\pi}{2N}}+\frac{1}{\sin^2\frac{N\pi}{2N}}\\
&\ \ \ +\frac{1}{\sin^2\frac{(N+1)\pi}{2N}}+\cdots+\frac{1}{\sin^2\frac{(2N-2)\pi}{2N}}+\frac{1}{\sin^2\frac{(2N-1)\pi}{2N}}\,.
\end{align*}
\end{small}
For $\sin x=\sin(\pi -x)$, $0 \leq x \leq 2\pi$ and $1/(\sin^2\frac{N\pi}{2N})=1$,
\begin{small}
\begin{align*}
\sum^{2N-1}_{n=1}\frac{1}{\sin^2\frac{n\pi}{2N}}=2\sum^{N-1}_{n=1}\frac{1}{\sin^2\frac{n\pi}{2N}}+1\,.
\end{align*}
\end{small}
Thus, we obtain identity (\ref{idty3}); that is,
\begin{equation*}
F_N(1)=\frac{1}{2N}\sum^{N-1}_{n=1}\frac{1}{\sin^2\frac{n\pi}{2N}}=\frac{1}{4N}\bigg(\frac{4N^2}{3}-\frac{1}{3}-1\bigg)=\frac{N}{3}-\frac{1}{3N}\,.
\end{equation*}

Similarly,
we expand (5) by noting that $\cos^2\frac{N\pi}{2N}=0$,
\begin{small}
\begin{align*}
\sum^{2N-1}_{n=1}\frac{\cos^2\frac{n\pi}{2N}}{\sin^2\frac{n\pi}{2N}}&=\sum^{N-1}_{n=1}\frac{\cos^2\frac{n\pi}{2N}}{\sin^2\frac{n\pi}{2N}}+\frac{\cos^2\frac{N\pi}{2N}}{\sin^2\frac{N\pi}{2N}}+\sum^{2N-1}_{n=N+1}\frac{\cos^2\frac{n\pi}{2N}}{\sin^2\frac{n\pi}{2N}}\\
&=2\sum^{N-1}_{n=1}\frac{\cos^2\frac{n\pi}{2N}}{\sin^2\frac{n\pi}{2N}}\,.
\end{align*}
\end{small}
For  $\sum^{N-1}_{n=1}\frac{\cos^2\phi_n}{\sin^2\phi_n}=
\frac{N^2}{3}-N+\frac{2}{3}$ , we have  $\sum^{2N-1}_{n=1}\frac{\cos^2\frac{n\pi}{2N}}{\sin^2\frac{n\pi}{2N}}=
\frac{4N^2}{3}-2N+\frac{2}{3}$. Therefore, we obtain identity (\ref{idty4}); that is,
\begin{equation*}
F_N(2)=\frac{2}{N}\sum^{N-1}_{n=1}\frac{\cos^2\phi_n/2}{\sin^2\phi_n/2}=
\frac{4N}{3}-2+\frac{2}{3N}\,.
\end{equation*}

By substituting (\ref{idty3}), (\ref{idty4}), and (\ref{retn2}) into (\ref{fSumForm}), we derive the following closed formula for $F_N(l)$
\begin{equation*}
F_N(l)=\frac{l^4}{12 N}-\frac{l^3}{3}+\frac{l^2 N}{3}-\frac{l^2}{6 N}+\frac{(-1)^l}{8
   N}+\frac{l}{3}-\frac{1}{8 N}.
\end{equation*}
By plugging $F_N(l)$ into (\ref{dfree}), we obtain the result in Proposition~\ref{path_bih_result}.
\end{IEEEproof}

\bibliographystyle{IEEEtran}
\bibliography{consensus}

\begin{thebibliography}{10}
\providecommand{\url}[1]{#1}
\csname url@rmstyle\endcsname
\providecommand{\newblock}{\relax}
\providecommand{\bibinfo}[2]{#2}
\providecommand\BIBentrySTDinterwordspacing{\spaceskip=0pt\relax}
\providecommand\BIBentryALTinterwordstretchfactor{4}
\providecommand\BIBentryALTinterwordspacing{\spaceskip=\fontdimen2\font plus
\BIBentryALTinterwordstretchfactor\fontdimen3\font minus
  \fontdimen4\font\relax}
\providecommand\BIBforeignlanguage[2]{{%
\expandafter\ifx\csname l@#1\endcsname\relax
\typeout{** WARNING: IEEEtran.bst: No hyphenation pattern has been}%
\typeout{** loaded for the language `#1'. Using the pattern for}%
\typeout{** the default language instead.}%
\else
\language=\csname l@#1\endcsname
\fi
#2}}

\bibitem{CaZa14}
R.~Carli and S.~Zampieri, ``Network clock synchronization based on the
  second-order linear consensus algorithm,'' \emph{IEEE Trans. Autom. Control},
  vol.~59, no.~2, pp. 409--422, 2014.

\bibitem{SuStBrGh15}
W.~Sun, E.~G. Str{\"o}m, F.~Br{\"a}nnstr{\"o}m, and M.~R. Gholami, ``Random
  broadcast based distributed consensus clock synchronization for mobile
  networks,'' \emph{IEEE Trans. Wireless Commun}, vol.~14, no.~6, pp.
  3378--3389, 2015.

\bibitem{DiFrMo99}
R.~Diekmann, A.~Frommer, and B.~Monien, ``Efficient schemes for nearest
  neighbor load balancing,'' \emph{Parallel Comput.}, vol.~25, no.~7, pp.
  789--812, 1999.

\bibitem{LiRu06}
Q.~Li and D.~Rus, ``Global clock synchronization in sensor networks,''
  \emph{IEEE Trans. Comput.}, vol.~55, no.~2, pp. 214--226, 2006.

\bibitem{FaMu04}
J.~A. Fax and R.~M. Murray, ``Information flow and cooperative control of
  vehicle formations,'' \emph{IEEE Trans. Autom. Control}, vol.~49, no.~9, pp.
  1465--1476, Sep. 2004.

\bibitem{Sa14}
A.~H. Sayed, ``Adaptation, learning, and optimization over networks,''
  \emph{Foundations and Trends® in Machine Learning}, vol.~7, no. 4-5, pp.
  311--801, 2014.

\bibitem{BaJoMiPa12}
B.~Bamieh, M.~R. Jovanovic, P.~Mitra, and S.~Patterson, ``Coherence in
  large-scale networks: Dimension-dependent limitations of local feedback,''
  \emph{{IEEE} Trans. Autom. Control}, vol.~57, no.~9, pp. 2235--2249, Sep.
  2012.

\bibitem{YoScLe10}
G.~F. Young, L.~Scardovi, and N.~E. Leonard, ``Robustness of noisy consensus
  dynamics with directed communication,'' in \emph{Proc. Amer. Control Conf.},
  Jun. 2010, pp. 6312--6317.

\bibitem{PaBa14}
S.~Patterson and B.~Bamieh, ``Consensus and coherence in fractal networks,''
  \emph{IEEE Trans. Control Netw. Syst.}, vol.~1, no.~4, pp. 338--348, Sep.
  2014.

\bibitem{FiLe16}
K.~Fitch and N.~E. Leonard, ``Joint centrality distinguishes optimal leaders in
  noisy networks,'' \emph{IEEE Trans. Control Netw. Syst.}, vol.~3, no.~4, pp.
  366--378, 2016.

\bibitem{YiZhLiCh15}
Y.~Yi, Z.~Zhang, Y.~Lin, and G.~Chen, ``Small-world topology can significantly
  improve the performance of noisy consensus in a complex network,''
  \emph{Comput. J.}, p. bxv014, 2015.

\bibitem{JaOl15}
\BIBentryALTinterwordspacing
A.~Jadbabaie and A.~Olshevsky, ``Scaling laws for consensus protocols subject
  to noise,'' \emph{arXiv:1508.00036}, 2015. [Online]. Available:
  \url{https://arxiv.org/abs/1508.00036}
\BIBentrySTDinterwordspacing

\bibitem{YoScLe16}
G.~F. Young, L.~Scardovi, and N.~E. Leonard, ``A new notion of effective
  resistance for directed graphs—part {I}: Definition and properties,''
  \emph{IEEE Trans. Autom. Control}, vol.~61, no.~7, pp. 1727--1736, 2016.

\bibitem{ReAt05}
W.~Ren and E.~Atkins, ``Second-order consensus protocols in multiple vehicle
  systems with local interactions,'' in \emph{AIAA Guidance, Navigation, and
  Control Conference and Exhibit}, 2005, pp. 15--18.

\bibitem{LiRuFu10}
Y.~Lipman, R.~M. Rustamov, and T.~A. Funkhouser, ``Biharmonic distance,''
  \emph{ACM Trans. Graph.}, vol.~29, no.~3, p.~27, 2010.

\bibitem{YiZhShCh17}
Y.~Yi, Z.~Zhang, L.~Shan, and G.~Chen, ``Robustness of first-and second-order
  consensus algorithms for a noisy scale-free small-world koch network,''
  \emph{IEEE Trans. Control Syst. Technol.}, vol.~25, no.~1, pp. 342--350,
  2017.

\bibitem{Pa17book}
G.~Patan{\`e}, ``An introduction to laplacian spectral distances and kernels:
  Theory, computation, and applications,'' \emph{Synthesis Lectures on Visual
  Computing: Computer Graphics, Animation, Computational Photography, and
  Imaging}, vol.~9, no.~2, pp. 1--139, 2017.

\bibitem{HuSzKo12}
D.~Hunt, B.~Szymanski, and G.~Korniss, ``Network coordination and
  synchronization in a noisy environment with time delays,'' \emph{Phys. Rev.
  E}, vol.~86, no.~5, p. 056114, 2012.

\bibitem{FiLe13}
K.~Fitch and N.~E. Leonard, ``Information centrality and optimal leader
  selection in noisy networks,'' in \emph{Proc. 52nd IEEE Conf. Decision
  Control}, 2013, pp. 7510--7515.

\bibitem{Xu09}
X.-P. Xu, ``Exact analytical results for quantum walks on star graphs,''
  \emph{J. Phys. A}, vol.~42, no.~11, p. 115205, 2009.

\bibitem{TzWu00}
W.-J. Tzeng and F.~Wu, ``Spanning trees on hypercubic lattices and
  nonorientable surfaces,'' \emph{Appl. Math. Lett.}, vol.~13, no.~7, pp.
  19--25, 2000.

\end{thebibliography}

\end{document}